\documentclass[11pt]{amsart}
\usepackage{amsmath,amsthm, amscd, amssymb, amsfonts, mathrsfs}
\usepackage[all]{xy}
\usepackage{color}
\usepackage{enumitem}
\usepackage{mathdots}%Please turn this on again after you finish

\usepackage{hhline} 
\usepackage[active]{srcltx}%{vista inversa}

\newcommand{\PSL}{\mathbf{PSL}}
\newcommand{\Gb}{\boldsymbol{G}}

\newcommand{\PSp}{\mathbf{PSp}}

\newcommand{\Pom}{\mathbf{P\Omega}}
\newcommand{\PSU}{\mathbf{PSU}}

\newcommand{\SU}{\mathbf{SU}}

\newcommand{\Jf}{{\tt J}}

\newcommand{\Ci}{\mathbf{C}}

\newcommand{\diag}{\operatorname{diag}}

\newcommand{\supp}{\operatorname{supp}}

\newcommand{\Inn}{\operatorname{Inn}}

\newcommand{\kc}{\mathbb F_q}

\newcommand\toba{\mathfrak B }
\newcommand\bpi{\boldsymbol{\pi}}

\newcommand{\trid}{\triangleright}

\newcommand{\Ha}{{\mathbb H}}

\newcommand{\kk}{\Bbbk}%{}

\newcommand{\Kb}{{\mathbb K}}

\newcommand{\Z}{\mathbb Z}
\newcommand{\N}{{\mathbb N}}
\newcommand{\I}{{\mathbb I}}

\newcommand{\G}{{\mathbb G}}
\newcommand{\B}{{\mathbb{B}}}
\newcommand{\T}{{\mathbb{T}}}
\newcommand{\U}{\mathbb{U}}
\newcommand\V{\mathbb V}
\newcommand{\Pa}{\mathbb{P}}

\newcommand{\M}{{\mathbb M}}

\newcommand{\F}{{\mathbb F}}

\newcommand{\GL}{\mathbf{GL}}
\newcommand{\GU}{\mathbf{GU}}

\newcommand{\SL}{\mathbf{SL}}
\newcommand{\Sp}{\mathbf{Sp}}
\newcommand{\SO}{\mathbf{SO}}

\newcommand{\Fr}{\operatorname{Fr}}

\newcommand{\Dc}{{\mathbf D}}

\newcommand{\Le}{{\mathbb L}}
\newcommand{\Vu}{{\mathbb V}}
\newcommand{\Oc}{{\mathcal O}}
\newcommand{\oc}{{\mathcal O}}

\newcommand\Gsc{\G_{\operatorname{sc}}}

\numberwithin{equation}{section}

\theoremstyle{plain}

\newtheorem{lema}{Lemma}[section]
\newtheorem*{maintheorem}{Main Theorem}
\newtheorem{theorem}[lema]{Theorem}

\newtheorem{prop}[lema]{Proposition}

\newtheorem{question-app}{Question}

\theoremstyle{definition}
\newtheorem{definition}[lema]{Definition}

\theoremstyle{remark}
\newtheorem{obs}[lema]{Remark}

\newcommand\id{\operatorname{id}}

\newcommand\s{\mathbb S}

\def\pf{\begin{proof}}
\def\epf{\end{proof}}

\theoremstyle{remark}

\newcounter{tabla}\stepcounter{tabla}

\begin{document}
%\[
\renewcommand{\baselinestretch}{1.2}

\thispagestyle{empty}
%\vspace*{2in}
\title[Nichols algebras over unipotent classes]
{Finite-dimensional pointed Hopf algebras\newline over finite simple groups of Lie type IV. \newline 
Unipotent classes in Chevalley and Steinberg groups}

\author[N. Andruskiewitsch, G. Carnovale, G. A. Garc\'ia]
{Nicol\'as Andruskiewitsch, Giovanna Carnovale and\newline Gast\'on Andr\'es Garc\'ia}

\thanks{
The work of N. A. was partially supported by CONICET, Secyt (UNC) and the
MathAmSud project GR2HOPF.	The work of G. C.
was partially supported by 
Progetto di Ateneo dell'Universit\`a di Padova: CPDA125818/12.
The work of G. A. G. 
was partially supported by CONICET, Secyt (UNLP) and ANPCyT-Foncyt 2014-0507.
The results were obtained during visits of N. A. and G. A. G. to the University of Padova, 
and of G. C. to the University of C\'ordoba, partially supported by the bilateral agreement between these Universities,
the Erasmus Mundus Action 2 Programme AMIDILA  and the Visiting Scientist program of the University of Padova.}

\address{\hspace{-20pt}  N. A. : FaMAF, 
Universidad Nacional de C\'ordoba. CIEM -- CONICET. %\newline \noindent
Medina Allende s/n (5000) Ciudad Universitaria, C\'ordoba,
Argentina}
\email{andrus@famaf.unc.edu.ar}
\address{\hspace{-20pt} G. C.:
Dipartimento di Matematica,
Universit\`a degli Studi di Padova,
via Trieste 63, 35121 Padova, Italia}
\email{carnoval@math.unipd.it}
\address{\hspace{-20pt} G. A. G.: Departamento de Matem\'atica, Facultad de Ciencias Exactas,
Universidad Nacional de La Plata. CONICET. C. C. 172, (1900)
La Plata, Argentina.}
\email{ggarcia@mate.unlp.edu.ar}

%\subjclass[2010]{16T05}
%\date{\today}

\begin{abstract}
We show that all unipotent classes in finite simple Chevalley or Steinberg groups, different from $\PSL_n(q)$ and
$\PSp_{2n}(q)$, collapse (i.e. are never the support of a finite-dimensional Nichols algebra), 
with a possible exception on one class of involutions in $\PSU_n(2^m)$.
\end{abstract}

\maketitle

\setcounter{tocdepth}{1}

\tableofcontents

\section{Introduction}

\subsection{The main result and the context} This is the fourth paper of our series on finite-dimensional complex
pointed Hopf algebras whose group of group-likes is isomorphic to a finite simple group of Lie type $\Gb$.
See Part I  \cite{ACG-I} for a comprehensive Introduction. As we explain in \emph{loc. cit.},
the primary task is to study Nichols algebras over $\Gb$ with support in a conjugacy class $\Oc$ of $\Gb$. 
Actually there are group-theoretical criteria allowing to conclude that every Nichols algebra with 
support in a given conjugacy class $\Oc$ has infinite dimension. These criteria were developed in \cite{AFGV-ampa,ACG-I,ACG-III}
and are recalled in \S \ref{subsec:racks}. The verification of any of these criteria in any conjugacy class might 
be difficult.
Let $p$ be a prime number, $m\in \N$, $q =p^m$, $\kc$ the field with $q$ elements and $\kk :=\overline{\kc}$.
There are three families of  finite simple groups of Lie type (according to the shape of the Steinberg endomorphism):
Chevalley, Steinberg and Suzuki-Ree groups;
see the list in \cite[p. 38]{ACG-I} and \cite[22.5]{MT} for  details.
Here are the contents of the previous papers:
\begin{itemize}[leftmargin=*]\renewcommand{\labelitemi}{$\diamond$}
\item In \cite{ACG-I} we dealt with  unipotent conjugacy classes in $\PSL_n(q)$, 
and as a consequence with the non-semisimple ones (since the 
centralizers of semisimple elements are products of groups with root system $A_{\ell}$).

\item The paper \cite{ACG-II} was devoted to  unipotent conjugacy classes in $\PSp_{2n}(q)$.

\item The subject of \cite{ACG-III} was the semisimple conjugacy classes in $\PSL_n(q)$. 
But we also introduced the criterium of type C, and applied it to some of the classes not reached with previous 
criteria in \cite{ACG-I,ACG-II}.
\end{itemize}

In this paper  we consider unipotent conjugacy classes in
Chevalley and Steinberg groups, different from $\PSL_n(q)$ and $\PSp_{2n}(q)$.  
Concretely, these are the groups in Table \ref{tab:groups}. Notice that $\PSU_3(2)$ is not simple
but needed for recursive arguments. 

\begin{table}[ht]
	\caption{Finite  groups considered in this paper; 
	$q$ odd for $\Pom_{2n+1}(q)$; $q\ge 3$ for $G_2(q)$}\label{tab:groups}
\begin{center}
\begin{tabular}{c|c c c|c}
\multicolumn{2}{c}{Chevalley} & & \multicolumn{2}{c}{Steinberg}
\\	
\cline{1-2}\cline{4-5}  $\Gb$  & \begin{small}Root system \end{small} &    & $\Gb$  & \begin{small}Root system \end{small}      \\
\cline{1-2}\cline{4-5}
$\Pom_{2n+1}(q)$ &  $B_{n}$, $n \ge 3$  
& &	$\PSU_n(q)$ &  $A_{n-1}$, $n \ge 3$    
\\
\cline{1-2}\cline{4-5}
$\Pom^+_{2n}(q)$& $D_{n}$, $n \ge 4$   
& &	$\Pom^-_{2n}(q)$& $D_{n}$, $n \ge 4$  
\\
\cline{1-2}\cline{4-5}
$G_2(q)$ & $G_{2}$   
& &	$^{3}D_4(q)$ & $D_{4}$  
\\
\cline{1-2}\cline{4-5}
$F_4(q)$ & $F_{4}$  
& &
	$^{2}E_6(q)$ & $E_{6}$ 
\\ 
\cline{1-2}\cline{4-5}
$E_j(q)$ & $E_{6}$, $E_7$, $E_8$ & & \multicolumn{2}{c}{ } \\
\cline{1-2}
\end{tabular}
\end{center}
\end{table}

\smallbreak  As in \cite[2.2]{AFGV-ampa}, we say that a conjugacy class $\oc$  of a finite group $G$
\emph{collapses} if the Nichols algebra $\toba(\Oc, \mathbf{q})$ has infinite dimension for
every finite faithful 2-cocycle $\mathbf{q}$. Our main result says:

\begin{maintheorem}
	Let $\Gb$ be  as in Table \ref{tab:groups}. Let $\Oc$ be a non-trivial unipotent conjugacy class in $\Gb$.
	Then either $\Oc$ collapses, or else $\Gb = \PSU_n(q)$ with $q$ even and $(2,1,\ldots,1)$  
	the partition of $\Oc$.
\end{maintheorem}

In the terminology of \S \ref{subsec:racks}, the classes not collapsing in the Main Theorem are
austere, see Lemma \ref{lem:SU211}. This means that the group-theoretical criteria do not apply for it; 
however, we ignore whether these classes collapse by other reasons. 
The classes in $\PSL_n(q)$ or $\PSp_{2n}(q)$ not collapsing (by these methods) are listed in Table \ref{tab:unip-slsp}.

\subsection{The scheme of the proof and organization of the paper}

Let $\Gb$ be a finite simple group of Lie type. Then  there is $q$ as above,
a simple simply connected algebraic group $\Gsc$ defined over  $\kc$
and a Steinberg endomorphism $F$ of $\Gsc$ such that $\Gb = \Gsc^F/Z(\Gsc^F)$. 
We refer to \cite[Chapter 21]{MT}
for details. Conversely, $\Gb = \Gsc^F/Z(\Gsc^F)$ is a simple group, out of a short list of exceptions, 
see \cite[Theorem 24.17]{MT}. 
For our inductive arguments, it is convenient to denote by $\Gb$ the quotient $\Gsc^F/Z(\Gsc^F)$ 
even when it is not simple.
Often there is a simple algebraic  group $\G$ with a projection $\bpi: \Gsc \to \G$
such that $F$ descends to $\G$ and  $[\G^F, \G^F]/ \bpi(Z(\Gsc^F)) \simeq \Gb$.

The proof of the Main Theorem is by application of the criteria of type C, D or F (see \S \ref{subsec:racks}), 
that hold by a recursive argument on the semisimple rank of $\Gsc$. 
The first step of the induction is given by the results on unipotent classes of $\PSL_n(q)$ and $\PSp_{2n}(q)$, 
while the recursive step is a reduction to Levi subgroups. 
Then we proceed group by group and class by class. The experience suggests that 
a general argument is not possible.
There are some exceptions in low rank for which Levi subgroups are too small and we 
need the representatives of the classes to apply \emph{ad-hoc} arguments. 

\smallbreak
Here is the organization of the paper:
We recall some notations and facts in \S \ref{sec:preliminaries}, where we also 
state the needed notation for groups of Lie type. In \S \ref{sec:unipotent} we describe the 
reduction to Levi subgroups and collect the known results on  unipotent classes of $\PSL_n(q)$ and $\PSp_{2n}(q)$. 

\smallbreak
Let $\Oc$ be a non-trivial unipotent class in a group $\Gb$ listed in Table \ref{tab:groups}. 
If $\Oc$ is not kthulhu then $\Oc$ collapses, cf. Theorem \ref{th:notkthulhu-colapses}.
The proof that $\Oc$ is not kthulhu is given in \S \ref{sec:chevalley},  respectively \S \ref{sec:steinberg},
when $\Gb$ is a  Chevalley, respectively Steinberg, group.

\smallbreak
Indeed, if 
$\Gb= \Pom_{2n+1}(q)$, $n\geq 3$, 
and $q$ odd, the claim is Proposition  \ref{prop:bn-qodd}. 
If $\Gb = \Pom^+_{2n}(q)$, $n \ge 4$,
$E_6(q)$, $E_7(q)$, or $E_8(q)$, then the claim is	 Proposition \ref{prop:DE}.
If $\Gb =F_4(q)$, the result follows from Lemmata \ref{lem:f4} and \ref{lem:f4-char2}; and 
if  $\Gb = G_2(q)$,  $q\geq 3$, the assertion follows from Lemmata \ref{lem:g2-qodd} and 
\ref{lem:g2-char2}. 

\smallbreak
In turn, $\PSU_n(q)$ is settled in  Proposition \ref{lem:sungral}; 
$\Pom^-_{2n}(q)$ in Proposition \ref{prop:bn-steinberg-gral}; 
$^{2}E_6(q)$ in Proposition \ref{prop:E6-steinberg}; 
and $^{3}D_4(q)$ in Proposition \ref{prop:triality-unipotent}.

\smallbreak
In this way, the Theorem is proved.

\subsection{Applications and perspectives}
 The results in this paper will be applied to settle the non-semisimple classes in Chevalley and Steinberg groups.
	
 Next we will deal with unipotent and non-semisimple classes in Suzuki-Ree groups. 
	These are too small to apply the recursive arguments introduced in this paper.
	
The semisimple conjugacy classes in $\Gb$ different from $\PSL_n(q)$ are more challenging. 
	We expect that reducible classes would collapse while the irreducible ones would be kthulhu, as is the case for 
	$\PSL_2(q)$ and $\PSL_3(q)$ (with some exceptions). Both cases require a deeper understanding of the classes,
	and in addition the irreducible case seems to need an inductive argument on the maximal subgroups.

\subsection*{Acknowledgements} At different stages of this project, Mauro Costantini, Paolo Papi and Jay Taylor, 
helped us with interesting conversations and precise references. We thank them a lot.

\section{Preliminaries}\label{sec:preliminaries}

If $a\leq b\in\N$, then $\I_{a,b}$ denotes $\{a,a+1,\,\ldots,b\}$;  for simplicity
$\I_a = \I_{1,a}$.

\subsection{Glossary of racks}\label{subsec:racks}
See \cite{ACG-III} for details and more information. 

\subsubsection{} A \emph{rack} is a  set $X \neq \emptyset$  with a self-distributive
operation $\trid: X \times X \to X$ such that $x\trid \underline{\quad}$ is
bijective for every $x \in X$. The archetypical example is the conjugacy class $\Oc^G_z$
of an element $z$ in a group $G$ with the operation $x\trid y = xyx^{-1}$, $x, y \in \Oc^G_z$. 
A rack $X$ is \emph{abelian} if $x\trid y = y$, for all $x, y \in X$.

\subsubsection{}  \cite[Definition 3.5]{AFGV-ampa} A rack $X$ is \emph{of type D}
if it has  a decomposable subrack
$Y = R\coprod S$ with elements $r\in R$, $s\in S$ such that
$r\trid(s\trid(r\trid s)) \neq s$. 

\begin{lema}\label{lema:ACG-2-10}\cite[Lemma 2.10]{ACG-I}
	Let $X$ and $Y$ be racks, $y_1\neq y_2\in Y$, $x_1\neq x_2\in X$
	such that  $x_1 \trid (x_2\trid (x_1\trid x_2)) \neq x_2$, $y_1\trid y_2 = y_2$.
	Then $X \times Y$ is of type D. \qed
\end{lema}

\begin{obs}\label{rem:real} One of 
	the hypothesis of Lemma \ref{lema:ACG-2-10} holds in the following setting.
Let $\Oc$ be a real conjugacy class, i.e. $\Oc=\Oc^{-1}$, with no involutions. 
Then $y_1 \neq y_2=y_1^{-1}$, that obviously commute. 
\end{obs}

\subsubsection{}  \cite[Definition 2.4]{ACG-I}
A rack $X$ is \emph{of type F}
if it has a family of subracks $(R_a)_{a \in \I_4}$
and elements $r_a\in R_a$, $a \in \I_4$, such that $R_a \triangleright R_b = R_b$,  for $a, b \in \I_4$, and
$R_a \cap R_b = \emptyset$, 
$r_a\triangleright r_b \neq r_b$ for $a\neq b \in \I_4$.

\subsubsection{}  \cite[Definition 2.3]{ACG-III}
A rack
$X$ is \emph{of type C}  when there are  a decomposable subrack
$Y = R\coprod S$  and elements $r\in R$, $s\in S$ such that $r \triangleright s \neq s$,
\begin{align*}
R &= \Oc^{\Inn Y}_{r},& S &= \Oc^{\Inn Y}_{s},&
\min \{\vert R \vert, \vert S \vert \}  &> 2 \text{ or } \max \{\vert R \vert, \vert S \vert \}  > 4.
\end{align*}
Here $\Inn Y$ is the subgroup of $\mathbb S_Y$ generated by $y\trid \underline{\quad}$, $y\in Y$.

\subsubsection{} 
Being of type C, D or F can be phrased in group terms, see \cite{ACG-III}. Here is a new formulation suitable for
later applications. 

\begin{lema}\label{lem:typeC}Let $\Oc$ be a conjugacy class in a group $H$. 
	If there are $r, s\in \Oc$ such that $r^2s\neq sr^2$, $s^2r\neq rs^2$ and 
	$\Oc_r^{\langle r,s\rangle}\neq\Oc_s^{\langle r,s\rangle}$, then 
	$\Oc$ is of type C.
\end{lema}
\pf We check that the conditions  in \cite[Lemma 2.8]{ACG-III} 
hold with $H=\langle r,s\rangle= \langle \Oc_r^{\langle r,s\rangle},\Oc_s^{\langle r,s\rangle}\rangle$. 
By hypothesis,  
$rs\neq sr$ and $\Oc_r^{\langle r,s\rangle}\neq\Oc_s^{\langle r,s\rangle}$. 
Now $r,s\trid r, s^2\trid r$ are all distinct, 
so $|\Oc_r^{\langle r,s\rangle}| >2$, and similarly for $\Oc_s^{\langle r,s\rangle}$.
\epf

\begin{theorem}\label{th:notkthulhu-colapses} \cite[Theorem  3.6]{AFGV-ampa},  \cite[Theorem 2.8]{ACG-I},  \cite[Theorem 2.9]{ACG-III}.
A rack $X$ of type D,  F or C collapses.
\end{theorem}

The proof rests on results from \cite{CH, HS, HV}.

\subsubsection{}\label{subsubsec:austere} A rack is
\begin{itemize}[leftmargin=*]\renewcommand{\labelitemi}{$\circ$}
	\item  \emph{kthulhu} if it is neither of type C, D nor F;
	
	\item  \emph{sober} if  every subrack is either abelian or indecomposable; 
		
	\item   \emph{austere} if every subrack generated by two elements is 
	either abelian or  indecomposable.
\end{itemize}
	
Clearly, sober implies austere and austere implies kthulhu. 

The criteria of type C, D, F are very flexible:
\begin{lema}\label{lem:flexibility} \cite{AFGV-ampa,ACG-I,ACG-III}
Let $Y$ be either a subrack or a quotient rack of a rack $X$. If $Y$ is not kthulhu, then $X$ is not kthulhu. \qed
\end{lema}

\subsection{Conjugacy classes}\label{subsec:conjugacy-classes}
\subsubsection{}\label{subsubsec:alg-gps}
Let $q = p^m$ be as above. We fix a simple algebraic group $\G$ defined over $\kc$, 
a maximal torus $\T$, with root system denoted by $\Phi$, and a Borel subgroup $\B$ containing $\T$.
We denote by $\U$ the unipotent radical of $\B$ and  by 
$\Delta \subset \Phi^+$ the corresponding sets of simple and positive roots.
Also  $\U^-$ is the unipotent radical of the opposite Borel subgroup $\B^-$ 
corresponding to $\Phi^-$.
We shall use the realisation of the associated root system and the numbering of simple roots  in \cite{Bou}. 
The coroot system of $\G$ is denoted by 
$\Phi^{\vee} = \{\beta^{\vee}~|~ \beta \in \Phi\} \subset X_{*}(\T)$,
where  $\langle\alpha, \beta^{\vee}\rangle = \frac{2(\alpha,\beta)}{(\beta,\beta)}$, 
for all $\alpha\in \Phi$.
Hence
\begin{align*}
\alpha(\beta^{\vee}(\zeta)) &= \zeta^{\frac{2(\alpha,\beta)}{(\beta,\beta)}}, 
& \alpha, \beta &\in \Phi, \zeta \in \kc^{\times}.
\end{align*}
We denote by $\Gsc$ the simply connected group covering $\G$.

For $\Pi\subset\Delta$,  we denote by $\Phi_{\Pi}$ the root subsystem with base $\Pi$ 
and $\Psi_\Pi:=\Phi^+ - \Phi_\Pi$.  For $\alpha\in\Phi$, we denote by 
$s_\alpha\in W=N_{\G}(\T)/\T$ the reflection with respect to $\alpha$. Also, $s_i = s_{\alpha_i}$, if $\alpha_i$ is a simple root with the alluded numeration.
Also, there is a
monomorphism of abelian groups $x_\alpha: \kk\to \U$; the image  $\U_{\alpha}$ of $x_\alpha$ 
is called a root subgroup.
We adopt the normalization of $x_\alpha$ and the notation for 
the elements in $\T$ from \cite[8.1.4]{springer}.
We recall the commutation rule:
$t\trid x_{\alpha}(a) = t x_{\alpha}(a)t^{-1}=x_{\alpha}(\alpha(t)a)$, for $t\in \T$ and $\alpha\in\Phi$. 
In particular, if $t= \beta^{\vee}(\xi)$ for some $\xi \in \Bbbk^{\times}$, then 
$t\trid x_{\alpha}(a) =x_{\alpha}(\alpha(\beta^{\vee}(\xi))a) =  x_{\alpha}(\xi^{\frac{2(\alpha,\beta)}{(\beta,\beta)}}a)$.

We denote by  $\Pa$ a standard parabolic subgroup of $\G$, with standard Levi subgroup $\Le$ and 
unipotent radical $\Vu$. 
Thus there exists $\Pi\subset\Delta$ such that $\Le = \langle \T, \U_{\pm\gamma}~|~ \gamma\in \Pi\rangle$. 

If $u\in\U$ then for every ordering of $\Phi^+$, there exist unique $c_\alpha\in \kk$ such that 
$u=\prod_{\alpha\in\Phi^+}x_{\alpha}(c_\alpha)$. We define
$\supp u=\{\alpha\in\Phi^+ ~|~ c_\alpha\neq0\}$. In general the 
support depends on the chosen ordering of $\Phi^+$. However, if $u\in \Vu$ as above, 
then $\supp u\subset \Psi_{\Pi}$ for every ordering of $\Phi^+$. 

\subsubsection{}\label{subsubsec:steinberg-endo} In this paper we deal with Chevalley and Steinberg groups.
Let $F$ be a Steinberg endomorphism of $\G$; it is the composition of the split endomorphism
$\Fr_q$ (the $q$-Frobenius map) with an automorphism induced by a Dynkin diagram automorphism $\vartheta$.
So, Chevalley groups correspond to $\vartheta = \id$.
We assume that $\T$ and $\B$ are F-stable.
Let $W^F = N_{\G^F}(\T) / \T^F$. 
Thus $W^F \simeq W$ for Chevalley groups.
For each $w \in W^F$, there is a representative $\dot w$ of $w$ in $N_{\G^F}(\T)$, cf. \cite[Proposition 23.2]{MT}.
Notice that   $\dot w \trid (\U_{\alpha}) = \U_{w(\alpha)}$ for all 
$\alpha \in \Phi$.
%\comgas{[We use later that $\supp (\dot w\trid u) = \{w(\alpha): \alpha \in \supp u\}$ for $w\in W$. It might be good that we state it here, do you agree?
%I understand that it is obvious, but I think it could help the non-expert.]}
Hence, if $\alpha, \beta \in \Phi$ are $\vartheta$-stable and have the same length, then
$\U_{\alpha}^F$ and $\U_{\beta}^F$ are conjugated by an element in $N_{\G^F}(\T)$ by \cite[Lemma 10.4 C]{HuLA}.
Let $\Gb = [\G^F, \G^F]/ Z(\G^F)$.

\subsubsection{}\label{alfabeta}
We shall often use the  Chevalley's commutator
formula \eqref{eq:chev}, see \cite[pp. 22 and  24]{yale}.
Let $\alpha, \beta \in\Phi$. If $\alpha+\beta$ is not a root, then  $\U_\alpha$ and $\U_\beta$ commute.
Assume that $\alpha+\beta\in\Phi$. Fix a total order in the set
$\Gamma$ of pairs $(i, j) \in \N ^2$ such that $i\alpha+j\beta\in \Phi$.
Then there exist $c_{ij}^{\alpha\beta}\in\kc$ such that
\begin{align}\label{eq:chev}x_{\alpha}(\xi)x_\beta(\eta)x_\alpha(\xi)^{-1}x_\beta(\eta)^{-1}
&= \prod_{(i,j)\in \Gamma}  x_{i\alpha + j\beta}(c_{ij}^{\alpha\beta}\xi^i\eta^j), &
\forall \xi, \eta&\in \kk.
\end{align}

\begin{definition}\label{def:alfabeta2} \cite[Definition 3.3]{ACG-II}
Let $\alpha,\beta \in \Phi^+$ such that $\alpha+\beta\in\Phi$ 
but the pair $\alpha,\beta$ does not appear in  Table  \ref{tab:c11-alpha-beta}.
We fix an ordering of $\Phi^+$.
A unipotent conjugacy class $\Oc$ in $\Gb$  has the
{\em $\alpha\beta$-property} if there exists $u\in \Oc\cap \U^F$ 
such that  $\alpha, \beta\in \supp u$ and for any expression 
$\alpha+\beta = \sum_{1\le i \le r} \gamma_i$, with  
$r > 1$ and $\gamma_i\in \supp u$, necessarily $r =2$ and $\{\gamma_1, \gamma_2\} = \{\alpha, \beta\}$.
\end{definition}

\begin{table}[ht]
	\caption{}\label{tab:c11-alpha-beta}
	\begin{align*}
&\begin{tabular}{|c|c|c|}
\multicolumn{3}{c}{$p = 3$}  \\
\hline    $\Phi$ & $\alpha$ & $\beta$  \\
\hline  \hline
 $G_2$ & $\alpha_1$ & $2\alpha_1+\alpha_2$ \\
\cline{2-3}
 & $2\alpha_1 +\alpha_2$ & $\alpha_1$ \\
\cline{2-3}
 & $\alpha_1 +\alpha_2$ & $2\alpha_1+\alpha_2$ \\
\cline{2-3}
 & $2\alpha_1+\alpha_2$ & $\alpha_1 +\alpha_2$ \\
\hline
\end{tabular}
&&\begin{tabular}{|c|c|c|}
\multicolumn{3}{c}{$p = 2$}  \\
	\hline  $\Phi$ & $\alpha$ & $\beta$  \\
	\hline  \hline
		 $B_n, C_n, F_4$  & \multicolumn{1}{r}{orthogonal} &to each other\\
	\hline
	 $G_2$ & $\alpha_1$ & $\alpha_1+\alpha_2$ \\
	\cline{2-3}
	& $\alpha_1+\alpha_2$ & $\alpha_1$ \\
	\hline
\end{tabular}
\end{align*}
\end{table}
Let $\alpha,\beta \in \Phi^+$.
The scalar $c_{1,1}^{\alpha,\beta}\neq0$ in \eqref{eq:chev}
if $\alpha+\beta\in\Phi$ and the pair does not appear in Table 
\ref{tab:c11-alpha-beta}.

\begin{prop}\label{prop:alfabetachevstein} \cite[Proposition 3.5]{ACG-II}
	Let $\Gb$ be a finite simple group of Lie type, with $q$ odd. Assume
	$\Oc$ has the $\alpha\beta$-property, for some $\alpha,\beta \in 
	\Phi^{+}$ such that $q>3$
	when $(\alpha,\beta)= 0$. Then  $\Oc$  is of type D. \qed
\end{prop}

\begin{obs}\label{obs:alfabeta-involution} 
	Assume $u$ satisfies the conditions in Definition \ref{def:alfabeta2}. 
	Then it is never an involution. Indeed if $q$ is 
	odd this is never the case. If $q$ is even then the 
	argument in the proof of \cite[Proposition 3.5]{ACG-II} shows that the coefficient
	of $x_{\alpha+\beta}$ in the expression of $u^2$ is nonzero. 
\end{obs}

\subsubsection{}\label{weyl-group-support} 
Let us choose an ordering of the positive roots and let $w\in W$ and $u\in \U$  
be such that $\Sigma:=w(\supp u)\subset\Phi^+$. 
Then, $\dot w\trid u\in \U$ and there is an ordering of the positive roots for 
which $\Sigma$ is the support of $\dot w\trid u$. 
If, in addition, $w\Sigma\not\subset\Psi_\Pi$ for some $\Pi\subset\Delta$, 
then by the discussion in \ref{subsubsec:alg-gps}, 
$\dot w\trid u\in \U - \Vu$.

\subsubsection{}\label{facts-root-systems} 
We shall need a fact on root systems.  Recall that there is a partial ordering 
$\preceq$ on the root lattice ${\mathbb Z}\Phi$ given by $\alpha\preceq\beta$ if 
$\beta-\alpha\in{\mathbb N}_0\Phi^+={\mathbb N}_0\Delta$. 
 
\begin{lema}\label{lem:rootsystems}Let $\gamma, \beta\in\Phi^+$ with $\beta\preceq \gamma$. 
	Then  there exists a sequence 
$\alpha_{i_1},\,\ldots,\alpha_{i_k}\in\Delta$ such that
 \begin{enumerate}[leftmargin=*]
  \item $\forall j\in \I_{k}$ we have $\gamma_j:=\beta+\alpha_{i_1}+\cdots+\alpha_{i_j}\in\Phi^+$;
  \item $\gamma=\gamma_k$.
 \end{enumerate}
 If, in addition, $\Phi$ is simply-laced, then $\gamma_j=s_{i_j}\cdots s_{i_1}\beta$ for every $j\in \I_{k}$. 
 \end{lema}
\pf (1) and (2) are consequences of \cite[Lemma 3.2]{sommers}, with $\alpha_1=\beta$, 
and the $\alpha_j$ being simple. 
Assume that $\Phi$ is simply-laced. Clearly, it is enough to prove it for a couple of roots. 
If $\alpha,\delta\in \Phi$ and $\alpha+\delta\in\Phi$, then
$\Phi \cap (\Z \alpha + \Z\delta)$ is a root system of type  $A_2$, so
$s_\alpha(\delta) = \alpha+\delta$. The last claim follows.
\epf

\section{Unipotent classes in finite groups of Lie type}\label{sec:unipotent}

\subsection{Reduction to Levi subgroups}

We start by Lemma \ref{lem:inductive-parabolic}, that  is behind the inductive step in most proofs below.
 We consider the following setting and notation, that we will use  throughout the paper:
\begin{itemize}[leftmargin=*]\renewcommand{\labelitemi}{}
	\item $\Pa_1,\ldots, \Pa_k$ are standard $F$-stable parabolic subgroups of $\G$;
	
	\item $\Pa_i=\Le_i\ltimes \Vu_i$ are Levi decompositions, with $\Le_i$ $F$-stable;
	
	\item $U_i:=(\U \cap \Le_i)^F$, $U_i^-:=(\U^- \cap \Le_i)^F$, $P_i:=\Pa_i^F$, $L_i:=\Le_i^F$, $V_i:=\Vu_i^F$; 
	
	\item $\pi_i: P_i\to L_i$ is the natural projection; 
	
	\item $M_i = \langle U_i, U_i^- \rangle \leq L_i$; for  $i \in \I_k$.
\end{itemize}

\begin{obs}\label{rem:def-Mi}
Assume that $\G=\G_{sc}$. If $\Le_i$ is standard, then $M_i = [\Le_i, \Le_i]^F$.
\end{obs}

\pf
Since $\G$ is simply-connected, so is $[\Le_i, \Le_i]$ (Borel-Tits, see \cite[Corollary 5.4]{sp-st}). 
Then \cite[Theorem 24.15]{MT} applies.
\epf 

%\comgas{[Somewhere we should recall \cite[Theorem 24.15]{MT}, that is
%$\G^{F} = \langle \U^{F}, (\U^{-})^{F}\rangle$, since we use it later, see for example
%the proof of Lemma 4.6]}

\begin{lema}\label{lem:inductive-parabolic} 
 
Let $u \in \U^F$. Then $\pi_i(u) \in M_i$ for all $i \in \I_k$.

\begin{enumerate}[leftmargin=*, label=\rm{(\alph*)}]
\item\label{item:parabolic} If $\Oc_{\pi_i(u)}^{M_i}$ 
is not kthulhu for some $i\in \I_k$, then $\Oc_{\pi_i(u)}^{L_i}$,  $\Oc_{u}^{P_i}$, and   
$\Oc_u^{\G^F}$ are not kthulhu either.
\item\label{item:all-parabolic} 
Assume that 
\begin{align}\label{eq:allccMi-notkthulhu}
\text{No non-trivial unipotent class in } M_i \text{ is kthulhu, } \forall i\in \I_k.
\end{align}
If $u\notin \cap_{i\in \I_k} V_i$, then $\Oc_{u}^{\G^F}$ is not kthulhu.
\item\label{item:all-parabolic-conjugate} 
Assume  that \eqref{eq:allccMi-notkthulhu} holds.  
Let $\Oc$ be a unipotent conjugacy class in $\G^F$.
If $\Oc\cap\U^F\not\subset\cap_{i\in \I_k} V_i$, then $\Oc$ is not kthulhu, hence collapses.
\end{enumerate}
\end{lema}
\pf Since $\U \le \Pa_i$, $u = u_1u_2$ with $u_1\in \Le_i$ and $u_2 \in \Vu_i \le \U$. Hence $u_1 \in \Le_i \cap \U$.
Since clearly $u_1$ and $u_2$ are $F$-invariant, $u_1 = \pi_i(u) \in M_i$.
Now 
\ref{item:parabolic} follows from Lemma \ref{lem:flexibility} 
and implies \ref{item:all-parabolic},  since $\pi_j(u)\not=1$ for some $j\in \I_k$. 
\ref{item:all-parabolic-conjugate} follows from \ref{item:all-parabolic} and Theorem \ref{th:notkthulhu-colapses}
because $\Oc\cap\U^F\not=\emptyset$.
\epf

\subsection{Unipotent classes in $\PSL_n(q)$ and $\PSp_{2n}(q)$}\label{subsec:psl-spu}

We recall now results in the previous papers of the series that constitute
the basis of the induction argument.
We will also need some of the non-simple groups of Lie 
type of small rank and small characteristic listed in \cite[3.2.1]{ACG-II}.

The following Theorem collects information from \cite[Table 2]{ACG-I}, \cite[Lemma 3.12 \& Tables 3, 4, 5]{ACG-II} 
and \cite[Tables 2 \& 3]{ACG-III}.

\begin{theorem}\label{thm:slsp}
	Let $\Gb$ be either $\PSL_n(q)$ or $\PSp_{2n}(q)$ and let $\Oc \neq \{e\}$ be a 
	unipotent conjugacy class in $\Gb$, not listed in Table \ref{tab:unip-slsp}. Then it is not kthulhu.
\end{theorem}
	\begin{table}[ht]
		\caption{Kthulhu classes in $\PSL_{n}(q)$ and $\PSp_{2n}(q)$}\label{tab:unip-slsp}
		\begin{center}
			\begin{tabular}{c|c|c}
				\hline $\Gb$   & class & $q$ \\
				\hline
				$\PSL_2(q)$   & $(2)$ & even, or $9$, or odd not a square   \\
				\hline
				$\PSL_3(2)$  &$(3)$  & $2$   \\
				\hline
				$\PSp_{2n}(q)$, $n\geq 2$  & $W(1)^{n-1}\oplus V(2)$  & even \\
				\hline
				$\PSp_{2n}(q)$, $n\geq2$  & $(2, 1^{2n-2})$  & $9$, or odd not a square\\
				\hline
				$\PSp_{4}(q)$  &  $W(2)$  & even  \\
				\hline
			\end{tabular}
		\end{center}
	\end{table}

We explain the notation of Table \ref{tab:unip-slsp}, see \cite{ACG-I, ACG-II} for further details:

\begin{enumerate}[leftmargin=*, label={(\roman*)}]
	\item Unipotent  classes in $\PSL_n(\kk)$ are parametrized by  partitions of $n$; i.e. $\lambda = (\lambda_1, \lambda_2, \dots)$ with $\lambda_1 \geq  \lambda_2 \geq \dots$ and $\sum_j \lambda_j = n$. Thus,
	$(n)$ is the regular unipotent class of $\PSL_n(\kk)$.
	Unipotent classes in $\PSL_n(q)$ with the same partition are isomorphic as racks.

	\item  Unipotent classes in 
	$\PSp_{2n}(\kk)$, for $q$ odd, are  also parametrized by  suitable partitions.

\item
Unipotent classes in $\PSp_{2n}(\kk)$, for $q$ even, are  
parametrized by their \emph{label}, which 
is the decomposition of the standard representation as a module for the action 
of an element in the conjugacy class:
\begin{align}\label{deco}
V &= \bigoplus_{i=1}^k W(m_i)^{a_i}\oplus \bigoplus_{j=1}^r V(2k_j)^{b_j},&  &0 <a_i,\; 0<b_j\leq 2,
\end{align}
for $m_i,k_j\geq1$. 
The block $W(m_i)$ corresponds to a unipotent class with partition $(m_i,m_i)$, 
whereas the block $V(2k_j)$ corresponds to a unipotent class with partition $(2k_j)$. 
\item The unipotent class in $\PSp_{4}(\kk)$ with label $W(2)$, respectively, in $\PSp_{2n}(\kk)$ 
with label $W(1)^{n-1}\oplus V(2)$ contain a unique unipotent class in $\PSp_{4}(q)$, respectively, $\PSp_{2n}(q)$. 
\end{enumerate}

\begin{obs}\label{obs:sp-not-involution}
Assume $q$ is even. If $\Oc$ is a unipotent conjugacy class in $\Sp_{2n}(q)$ enjoying the  
$\alpha\beta$-property, for some $\alpha$ and $\beta$,  then $\Oc$ is of type C, D, or F. 
Indeed, by Theorem  \ref{thm:slsp} the kthulhu unipotent classes in $\Sp_{2n}(q)$ for $q$ 
even consists of involutions.  Remark \ref{obs:alfabeta-involution} applies. 
\end{obs}

\subsection{Further remarks}\label{subsec:further-rem}
If a product $X = X_1 \times X_2$  of racks has a factor $X_1$ that is not kthulhu, then neither is $X$.
Indeed, pick $x \in X_2$ ; then $X_1 \times \{x\}$ is a subrack of $X$ and Lemma \ref{lem:flexibility} applies
(here as usual $X_2$ can be realized as a subrack of a group,  so that $x \trid x = x$).
The following results will be needed in order to deal with products of possibly kthulhu racks.

\begin{lema}\label{lem:x1x2y1y2}
Let $\Oc$ be a unipotent conjugacy class in Table \ref{tab:unip-slsp}.
\begin{enumerate}[leftmargin=*, label={\rm(\alph*)}]
	\item\label{item:x1x2y1y2-a} There exist $x_1$, $x_2\in \Oc$ such that $(x_1x_2)^2\neq (x_2x_1)^2$. 
	
	\item\label{item:x1x2y1y2-b}
	If $\Gb \neq \PSL_2(2),\PSL_2(3)$, then 
	there exist $y_1,y_2\in \Oc$ such that $y_{1}\neq y_{2}$ and $y_1y_2= y_2y_1$.
\end{enumerate}
\end{lema}

\pf By the isogeny argument \cite[Lemma 1.2]{ACG-I}, we may reduce to classes 
in $\SL_n(q)$ or $\Sp_{2n}(q)$. Also,
the classes in $\PSp_{4}(q)$ with label $W(2)$ and $W(1)\oplus V(2)$ are isomorphic 
as racks, \cite[Lemma 4.26]{ACG-II}, so we need not to deal with the last row in Table  \ref{tab:unip-slsp}.

If $\Oc$ is the class in $\SL_3(2)$, then $x_1=\id+e_{1,2}+e_{2,3}$  and $x_2=\sigma\trid x_1$, where
$\sigma=e_{1,2}+e_{2,1}+e_{3,3}$, do the job for  \ref{item:x1x2y1y2-a}. For \ref{item:x1x2y1y2-b},
take $y_1 = x_1$ and $y_2=x_1^3=x_1^{-1}$, that belongs to $\Oc$ by  \cite[Lemma 3.3]{ACG-I}.

%Assume that $\Oc$ is the class with label $W(2)$ in $\Sp_4(2)\simeq \s_6$. Then 
%$x_1=\id+e_{1,2}+e_{3,4}$ and $x_2=\sigma\trid x_1$, where $\sigma:=e_{1,2}+e_{2,1}+e_{3,4}+e_{4,3}$,
%are as needed in \ref{item:x1x2y1y2-a}. Also  all classes of involutions in $\s_6$ 
%have the property \ref{item:x1x2y1y2-b}.

If $\Oc$ is the class in $\SL_2(q)$, then $x_1=\id+e_{1,2} \in \Oc$ and  $x_2=\sigma\trid x_1$, where
$\sigma=e_{1,2}-e_{2,1}$ do the job for  \ref{item:x1x2y1y2-a}; 
while $y_1=x_1$, and $y_2=\id+a^2e_{1,2}$,  for $a\in \kc$, $a^2\neq0,1$, are as needed in \ref{item:x1x2y1y2-b}
when $q > 3$.

Finally, let $\Oc$ be one of the classes in $\Sp_{2n}(q)$, cf. Table \ref{tab:unip-slsp}. 
Then $x_1=\id+e_{1,2n} \in \Oc$ and  $x_2=\sigma\trid x_1$, where
$\sigma=e_{1,2n}-e_{2n,1}+\sum_{j\neq 1,2n}e_{jj}$ do the job for  \ref{item:x1x2y1y2-a}.
Let $\tau$ be the block-diagonal matrix
$\tau=\diag(\Jf_2,\id_{2n-2},\Jf_2)$, with $\Jf_2=\left(\begin{smallmatrix}0&1\\
1&0\end{smallmatrix}\right)$. Then  $\tau\in \Sp_{2n}(q)$ and $y_1:=x_1$, $y_2:=\tau\trid y_1$ fulfil \ref{item:x1x2y1y2-b}.
\epf

Here are results on regular unipotent  classes needed later.
Let $\Gsc$ be a simply connected simple algebraic group and $F$ a Steinberg endomorphism as before;
let $\Gb = \Gsc^F/ Z(\Gsc^F)$ but we do not assume that $\Gb$ is simple.

\begin{prop}\label{prop:regular}\cite[3.7, 3.8, 3.11]{ACG-II} Let $\Oc$ be a regular 
unipotent  class in $\Gb$. 
If any of the conditions below is satisfied, then $\Oc$ is of type  D, or F.
\begin{enumerate}[leftmargin=*]
\item\label{item:regular-chevalley} $\Gb\neq\PSL_2(q)$ is Chevalley and $q\neq 2,4$;
\item\label{item:regular-psu3} $\Gb=\PSU_3(q)$, with $q\neq 2,8$;  
\item\label{item:regular-psu4} $\Gb=\PSU_4(q)$, with $q\neq 2,4$;  
\item\label{item:regular-psun} $\Gb=\PSU_n(q)$, with $n\geq5$ and $q\neq2$;
\end{enumerate}
In addition, every regular unipotent class in $\GU_n(q)$,  where $1<n$ is odd and $q=2^{2h+1}$,  $h\in \N_0$, is of type D. \qed
\end{prop}
 
Finally, we quote \cite[Lemma 4.8]{ACG-II}:
 
 \begin{lema}\label{lem:x1x2}Let $\Oc$ be a regular unipotent class in 
either $\SL_n(q)$, $\SU_n(q)$ or $\Sp_{2n}(q)$, $q$ even. 
Then there are $x_1,x_2\in\Oc$ such that $(x_1x_2)^2\neq(x_2x_1)^2$.
\end{lema}

\section{Unipotent classes in Chevalley groups}\label{sec:chevalley}

In this Section we deal with unipotent conjugacy classes in a finite simple Chevalley group $\Gb = \Gsc^F/Z(\Gsc^F)$, different 
to $\PSL_n(q)$ and $\PSp_{2n}(q)$, treated in \cite{ACG-I,ACG-II}, see \S \ref{subsec:psl-spu}. 
For convenience, we shall work in $\G_{sc}^{F}$, cf. \cite[Lemma 1.2]{ACG-I}.
Let 
\begin{align}
\Psi(\beta) &=\{\gamma\in\Phi~|~\beta\preceq\gamma\}, & \beta&\in\Phi.
\end{align}

Let $u\in \U$ and $ \beta \in\Phi^+$. Then the support $\supp u$ depends on a  fixed ordering of $\Phi^+$, but the assertion 
$\supp u \subset \Psi(\beta)$ does not. Indeed, passing from one order to another boils down to 
successive applications of the Chevalley formula \eqref{eq:chev}, that do not affect the claim.

\smallbreak
We denote by $\Oc$ a non-trivial unipotent conjugacy class in $\Gb$.

\subsection{Unipotent classes in $\Pom^+_{2n}(q)$, $n \ge 4$; $E_6(q)$, $E_7(q)$ 
and $E_8(q)$}\label{sec:unipotent-notFG}
We first deal with the case when $\Phi$ simply-laced, i.e. $\Gb$ is one of
 $\Pom^+_{2n}(q)$, $n \ge 4$; $E_6(q)$, $E_7(q)$ and $E_8(q)$.

\begin{lema}\label{lem:supports}  Given $\beta\in\Phi^+ - \Delta$, there is $x\in \Oc\cap\U^F$ 
with $\supp x\not\subset\Psi(\beta)$. \end{lema}

\pf  Let $u\in\Oc\cap\U^F$. 
If $\supp u\not\subset\Psi(\beta)$, then we are done. Assume that $\supp u\subset\Psi(\beta)$. 
We claim that there is $\tau\in N_{\G_{sc}^F}(\T)$ such that 
\begin{align*}
x &:=\tau\trid u\in \Oc\cap \U^F& &\text{ and } &&\supp x\not\subset\Psi(\beta).
\end{align*}
For every $\gamma\in\Psi(\beta)$ there is a unique $k$ such that $\gamma=\beta+\alpha_{i_1}+\cdots+\alpha_{i_k}$ 
as in Lemma \ref{lem:rootsystems}. Let $m$ be the minimum $k$ for $\gamma\in \supp u$.
We call $m$ the bound of $u$. We will prove the claim by induction on the bound $m$.   
If $m=0$ then $\beta\in \supp u$ and since $\beta\not\in\Delta$, there is a simple reflection $s_i$ 
such that $s_i\beta=\beta-\alpha_i\in\Phi^+ - \Psi(\beta)$. Also, $s_i\gamma\in\Phi^+$ 
for every $\gamma\in \supp u$ because $s_i(\Phi^+ - \{\alpha_i\})=\Phi^+ - \{\alpha_i\}$. 
In this case we take $\tau= \dot {s_i}$ to be any representative of $s_i$ in 
$N_{\G_{sc}^{F}}(\T)$.

Let now $m>0$ and assume that the statement is proved for unipotent elements with bound $m-1$. 
Let $\gamma\in\supp u$ reach the minimum, i.e., be such that $\gamma=\beta+\alpha_{i_1}+\cdots+\alpha_{i_{m}}$ 
for some $\alpha_{i_j}\in\Delta$ chosen as in Lemma \ref{lem:rootsystems}. Then  
$\gamma'=s_{i_{m}}\gamma=\beta+\alpha_{i_1}+\cdots+\alpha_{i_{m-1}}\in\Psi(\beta)$, and 
$s_{i_m}\alpha\in\Phi^+$ for every $\alpha\in\Psi(\beta)$ by construction. Let $\dot{s}_{i_m}$ 
be a representative of $s_{i_{m}}$ in $N_{\Gsc^F}(\T)$. Then 
$u'=\dot{s}_{i_m}\trid u\in\Oc\cap \U^F$ and either $\supp u'\not\subset\Psi(\beta)$, or $\supp u'\subset\Psi(\beta)$, 
with bound at most $m-1$. 
In the first case, we conclude by setting $x=u'$. In the second case, we use the inductive hypothesis.\epf

\begin{prop}\label{prop:DE}  $\Oc$ is not kthulhu. 
\end{prop}
\pf  The basic idea of the proof is to apply  
Lemma \ref{lem:inductive-parabolic} \ref{item:all-parabolic-conjugate}
to a series of standard $F$-stable parabolic subgroups $\Pa_i$ of $\Gsc$ for 
which \eqref{eq:allccMi-notkthulhu} holds. We show that for every $\Oc$ and for every $\Gb$,  
we have $\Oc\cap \U^F\not\subset \cap_iV_i$. 
This follows from Lemma \ref{lem:supports} by observing that in each case $\cap_iV_i$ is a product of root subgroups 
corresponding to roots in $\Psi(\beta)$ for 
some $\beta\in\Phi^+ - \Delta$. 
We analyze the different cases according to  $\Phi$.

\subsubsection*{$D_n$, $n\geq4$}  We consider the parabolic subgroups $\Pa_1$ and $\Pa_2$ such 
that $\Le_1$ and $\Le_2$ have root systems $A_{n-1}$, 
generated respectively by $\Delta -  \alpha_{n-1}$  and  $\Delta -  \alpha_{n}$. Since $n\geq 4$, \eqref{eq:allccMi-notkthulhu} 
holds by Theorem \ref{thm:slsp}. Let $u \in V_1\cap V_2$. Then $\alpha \in \supp u$ 
if and only if $\alpha$ contains $\alpha_{n-1}$ and $\alpha_n$ in its expression, i.e.
$\alpha\in \Psi(\beta)$ for $\beta=\alpha_{n-2}+\alpha_{n-1}+\alpha_n$. By Lemma \ref{lem:supports}, 
$\Oc\cap \U^F\not\subset \cap_iV_i$.

\subsubsection*{$E_6$} We consider the parabolic subgroups $\Pa_1$, $\Pa_2$ and 
$\Pa_3$ such that $\Le_1$, $\Le_2$ and $\Le_3$ have root systems $D_5$, $D_5$ and $A_5$,
generated respectively by $\Delta -  \alpha_1$,    $\Delta -  \alpha_6$  and  $\Delta -  \alpha_{2}$. 
By Theorem \ref{thm:slsp} and the result for $D_n$, \eqref{eq:allccMi-notkthulhu}  holds. 
Let $u \in V_1\cap V_2 \cap V_3$. Then $\alpha \in \supp u$ 
if and only if 
$\alpha\in \Psi(\beta)$ for $\beta=\sum_{i=1}^6\alpha_i$. By Lemma \ref{lem:supports}, 
$\Oc\cap \U^F\not\subset \cap_iV_i$.

\subsubsection*{$E_7$} We consider the parabolic subgroups $\Pa_1$, $\Pa_2$ and $\Pa_3$ 
such that $\Le_1$, $\Le_2$ and $\Le_3$ have root systems $D_6$, $E_6$ and $A_6$,
generated respectively by $\Delta -  \alpha_1$,    $\Delta -  \alpha_7$  and  $\Delta -  \alpha_{2}$. 
By Theorem \ref{thm:slsp} and the results for $D_n$ and $E_6$, \eqref{eq:allccMi-notkthulhu}  holds. 
Let $u \in V_1\cap V_2 \cap V_3$. Then $\alpha \in \supp u$ 
if and only if 
$\alpha\in \Psi(\beta)$ for $\beta=\sum_{i=1}^7\alpha_i$. By Lemma \ref{lem:supports}, 
$\Oc\cap \U^F\not\subset \cap_iV_i$.

\subsubsection*{$E_8$} We consider the parabolic subgroups $\Pa_1$, $\Pa_2$ and 
$\Pa_3$ 
such that $\Le_1$, $\Le_2$ and $\Le_3$ have root systems $D_7$, $E_7$ and $A_7$,
generated respectively by $\Delta -  \alpha_1$,    $\Delta -  \alpha_8$  and  $\Delta -  \alpha_{2}$. 
By Theorem \ref{thm:slsp} and the results for $D_n$ and $E_7$, \eqref{eq:allccMi-notkthulhu}  holds. 
Let $u \in V_1\cap V_2 \cap V_3$. Then $\alpha \in \supp u$ 
if and only if 
$\alpha\in \Psi(\beta)$ for $\beta=\sum_{i=1}^8\alpha_i$. By Lemma \ref{lem:supports}, 
$\Oc\cap \U^F\not\subset \cap_iV_i$. 
\epf

\subsection{Unipotent classes in $\Pom_{2n+1}(q)$}
Here we  deal with  $\Pom_{2n+1}(q)$, i.e. $\Phi$ is of type $B_n$,  $n\geq 3$. 
In this case, $q$ is always odd.

\begin{prop}\label{prop:bn-qodd} $\Oc$ is not kthulhu.  
\end{prop}
\pf We consider the standard $F$-stable parabolic subgroups $\Pa_1$ and $\Pa_2$ such 
that $\Le_1$ and $\Le_2$ have root systems $A_{n-1}$ and $C_2$, 
generated respectively by $\Pi_1:=\Delta -  \alpha_{n}$  and  $\Pi_2=\{\alpha_{n-1},\alpha_n\}$. 
By Lemma \ref{lem:inductive-parabolic} \ref{item:parabolic} and Theorem \ref{thm:slsp}, 
if $\Oc\cap \U^F\not\subset V_1$  then $\Oc$ is not kthulhu.
Let us thus consider $u\in\Oc\cap V_1$. Then $\supp u\subset \Psi_{\Pi_1}= 
\{\varepsilon_i, \varepsilon_j+\varepsilon_l~|~ i,j,l\in\I_n,  j<l\}$, since it must contain
$\alpha_{n}$. We will apply the argument in \ref{weyl-group-support}.

Assume first that $\supp u\subset \{\varepsilon_j+\varepsilon_l~|~ j,l\in\I_n,  j<l\}$. 
Let  $\ell$ be the maximum ${l}$ such that $\varepsilon_j+\varepsilon_l\in \supp u$ for some $j\in\I_{n-1}$.
Then $s_{\varepsilon_\ell}(\supp u)\subset\Phi^+$. Let $\dot{s}_{\varepsilon_\ell}$ 
be a representative of $s_{\varepsilon_\ell}$ in $N_{\Gsc^F}(\T)$.  Then 
$\dot{s}_{\varepsilon_\ell}\trid u\in\Oc\cap \U^F$ and 
$\varepsilon_j-\varepsilon_\ell\in \supp(\dot{s}_{\varepsilon_\ell}\trid u)$ for every $j$ such that 
$\varepsilon_j+\varepsilon_\ell\in \supp u$. 
Hence $\dot{s}_{\varepsilon_{\ell}}\trid u\in\Oc \cap \U^F-  V_1$. 
By the previous argument, $\Oc$ is not kthulhu.

Assume next that there is some $i$ such that $\varepsilon_i\in \supp u$. 
We can always assume $i=n$. Indeed, if $\varepsilon_n\not\in\supp u$, 
we may replace $u$ by $\dot{s}_{\varepsilon_i-\varepsilon_n}\trid u\in\Oc\cap \U^F$,
where $\dot{s}_{\varepsilon_i-\varepsilon_n}$ is a representative of $s_{\varepsilon_i-\varepsilon_n}$ in 
$N_{\Gsc^F}(\T)$. Then
$\pi_{2}(u)\in M_2$ lies in a non-trivial unipotent conjugacy class in a group isomorphic to $\Sp_{4}(q)$ 
and the short simple root lies in the support. 
A direct computation shows that a representative of this class in $\Sp_4(q)$ is as follows:
\begin{align*}
&\left(\begin{smallmatrix}
1&a&*&*\\
0&1&0&*\\
0&0&1&-a\\
0&0&0&1
\end{smallmatrix}\right),&  a&\neq 0.
\end{align*} 
Thus, its Jordan form has partition $(2,2)$ and this class is not kthulhu by  
Theorem \ref{thm:slsp} (recall that $q$ is odd). Then Lemma \ref{lem:inductive-parabolic} applies.
\epf

\subsection{Unipotent classes in  $F_4(q)$}\label{subsec:f4}
\begin{table}[ht]
	\caption{Representatives of unipotent classes in $F_4(q)$ in odd characteristic; $\eta,\xi$ and $\zeta$ are suitable elements in $\kc^\times$}\label{tab:f4-odd-char-shoji}
	\begin{center}
		\begin{tabular}{|l|}
			\hline
			$x_1=x_{1+2}(1)$\\
			$x_2=x_{1-2}(1)x_{1+2}(-1)$\\
			$x_3=x_{1-2}(1)x_{1+2}(-\eta)$\\
			$x_4=x_{2}(1)x_{3+4}(1)$\\
			$x_5=x_{2-3}(1)x_{4}(1)x_{2+3}(1)$\\
			$x_6=x_{2-3}(1)x_{4}(1)x_{2+3}(\eta)$\\
			$x_7=x_{2}(1)x_{1-2+3+4}(1)$ \\
			$x_8=x_{2-3}(1)x_{4}(1)x_{1-2}(1)$\\
			$x_9=x_{2-3}(1)x_{3-4}(1)x_{3+4}(-1)$ \\
			$x_{10}=x_{2-3}(1)x_{3-4}(1)x_{3+4}(-\eta)$\\
			$x_{11}=x_{2+3}(1)x_{1+2-3-4}(1)x_{1-2+3+4}(1)$\\
			$x_{12}=x_{2-3}(1)x_{4}(1)x_{1-4}(1)$\\
			$x_{13}=x_{2-3}(1)x_{4}(1)x_{1-4}(\eta)$\\
			$x_{14}=x_{2-4}(1)x_{3+4}(1)x_{1-2}(-1)x_{1-3}(-1)$\\
			$x_{15}=x_{2-4}(1)x_{3+4}(1)x_{1-2}(-\eta)x_{1-3}(-1)$ \\
			$x_{16}=x_{2-4}(1)x_{2+4}(-\eta)x_{1-2+3+4}(1)x_{1-3}(-1)$\\
			$x_{17}=x_{2-4}(1)x_{3+4}(1)x_{1-2-3+4}(1)x_{1-2}(-\eta)x_{1-3}(\xi)$\\
			$x_{18}=x_{2}(1)x_{3+4}(1)x_{1-2+3-4}(1)x_{1-2}(-1)x_{1-3}(\zeta)$\\
			$x_{19}=x_{2-3}(1)x_{3-4}(1)x_{4}(1)$\\
			$x_{20}=x_{2}(1)x_{3+4}(1)x_{1-2-3-4}(1)$ \\
			$x_{21}=x_{2-4}(1)x_{3}(1)x_{2+4}(1)x_{1-2-3+4}(1)$ \\
			$x_{22}=x_{2-4}(1)x_{3}(1)x_{2+4}(\eta)x_{1-2-3+4}(1)$ \\
			$x_{23}=x_{2-3}(1)x_{3-4}(1)x_{4}(1)x_{1-2}(1)$\\
			$x_{24}=x_{2-3}(1)x_{3-4}(1)x_{4}(1)x_{1-2}(\eta)$\\
			$x_{25}=x_{2-3}(1)x_{3-4}(1)x_{4}(1)x_{1-2-3-4}(1)$\\
			$x_{26}=x_{2-3}(1)x_{3-4}(1)x_{4}(1)x_{1-2-3-4}(1)x_{1-2+3+4}(\zeta)$\\
			$x_{27}=x_{2-3}(1)x_{3-4}(1)x_{4}(1)x_{1-2-3-4}(1)x_{1-2+3+4}(-\zeta)$\\
			\hline
		\end{tabular}
	\end{center}
\end{table}

\begin{table}[ht]
	\caption{Representatives of unipotent classes in $F_4(q)$ in even characteristic; $\eta$ and $\zeta$ are suitable elements in $\kc^\times$}\label{tab:f4-char2-shinoda}
	\begin{center}
		\begin{tabular}{|l|}
			\hline
			$x_1=x_{1}(1)$\\
			$x_2=x_{1+2}(1)$\\
			$x_3=x_{1}(1)x_{1+2}(1)$\\
			$x_4=x_{2+3}(1)x_{1}(1)$\\
			$x_5=x_{2}(1)x_{2+3}(1)x_{1-3}(1)$\\
			$x_6=x_{2}(1)x_{2+3}(1)x_{1-3}(1)x_{1+3}(\eta)$\\
			$x_7=x_{2}(1) x_{2+3}(1)x_{1-2+3+4}(1)$ \\
			$x_8=x_{2}(1) x_{2+3}(1)x_{1-2+3+4}(1)x_{1+4}(\eta)$ \\
			$x_9=x_{2}(1)x_{1-2}(1)$\\
			$x_{10}=x_{2}(1)x_{1-2}(1)x_{1+2}(\eta)$\\
			$x_{11}=x_{2}(1)x_{3+4}(1)x_{1-4}(1)$\\
			$x_{12}=x_{2}(1)x_{1-2+3+4}(1)x_{1-4}(1)$\\
			$x_{13}=x_{2}(1)x_{2+3}(1)x_{1-2}(1)$\\
			$x_{14}=x_{2}(1)x_{3+4}(1)x_{1-2}(1)$\\
			$x_{15}=x_{2}(1)x_{2+3}(1)x_{1-2+3+4}(1)x_{1-3}(1)$\\
			$x_{16}=x_{2}(1)x_{2+3}(1)x_{1-2+3+4}(1)x_{1-2}(1)$\\
			$x_{17}=x_{2}(1)x_{2+3}(1)x_{1-2-3+4}(1)x_{1-2}(1)$\\
			$x_{18}=x_{2}(1)x_{2+3}(1)x_{1-2-3+4}(1)x_{1-2}(1)x_{1-4}(\eta)$\\
			$x_{19}=x_{2}(1)x_{3+4}(1)x_{1-2+3-4}(1)x_{1-2}(1)x_{1-3}(\zeta)$\\
			$x_{20}=x_{1-2}(1)x_{2-3}(1)x_3(1)$\\
			$x_{21}=x_{1-2}(1)x_{2-3}(1)x_3(1)x_{2+3}(\eta)$\\
			$x_{22}=x_{4}(1)x_{2-4}(1) x_{1-2+3-4}(1)$\\
			$x_{23}=x_{4}(1)x_{2-4}(1) x_{2+4}(\eta)x_{1-2+3-4}(1)$\\
			$x_{24}=x_{2-4}(1)x_{3+4}(1)x_{1-2-3-4}(1)x_{1-2-3+4}(1)$\\
			$x_{25}=x_{2-4}(1)x_{3+4}(1)x_{1-2-3-4}(1)x_{1-2-3+4}(1)x_{1-2}(\eta)$\\
			$x_{26}=x_{2-4}(1)x_{3+4}(1)x_{1-2-3-4}(1)x_{1-2-3+4}(1)x_{1-2}(\eta)x_{1-3}(\eta)$\\
			$x_{27}=x_{2-4}(1)x_{3}(1)x_{3+4}(1)x_{2+4}(\eta)x_{1-2-3+4}(1)$\\
			$x_{28}=x_{2-4}(1)x_{3}(1)x_{3+4}(1)x_{2+4}(\eta)x_{1-2-3+4}(1)x_{1-2}(\eta)$\\
			$x_{29}=x_{1-2}(1)x_{2-3}(1)x_{3-4}(1)x_{4}(1)$\\
			$x_{30}=x_{1-2}(1)x_{2-3}(1)x_{3-4}(1)x_{4}(1)x_{3+4}(\eta)$\\
			$x_{31}=x_{2-3}(1)x_{3-4}(1)x_{4}(1)x_{1-2-3-4}(1)$\\
			$x_{32}=x_{2-3}(1)x_{3-4}(1)x_{4}(1)x_{1-2-3-4}(1)x_{3+4}(\eta)$\\
			$x_{33}=x_{2-3}(1)x_{3-4}(1)x_{4}(1)x_{1-2-3-4}(1)x_{1-2}(\eta)$\\
			$x_{34}=x_{2-3}(1)x_{3-4}(1)x_{4}(1)x_{1-2-3-4}(1)x_{3+4}(\eta)x_{1-2}(\eta)$\\
			\hline
		\end{tabular}
	\end{center}
\end{table}

Here we deal with unipotent classes in  $F_4(q)$. In this case the approach in 
Section \ref{sec:unipotent-notFG} is not effective.
Indeed, in characteristic $2$, \eqref{eq:allccMi-notkthulhu} does not hold for any of the standard parabolic subgroups. 
For this reason we shall use 
explicit representatives of unipotent classes and apply results from Theorem \ref{thm:slsp} and Proposition  
\ref{prop:bn-qodd} for $B_3$, where $q$ is assume to be odd.

We use the list of representatives of unipotent classes in $F_4(q)$ in \cite[Tables 5,6]{shoji} for $q$ odd, see Table \ref{tab:f4-odd-char-shoji}, respectively in \cite[Theorem 2.1]{shinoda} for $q$ even, see Table \ref{tab:f4-char2-shinoda}. 
We indicate the roots  as in \cite{shinoda}:
$\varepsilon_i$ is indicated by $i$, $\varepsilon_i-\varepsilon_j$ 
is indicated by $i-j$, and $\frac{1}{2}(\varepsilon_1\pm\varepsilon_2\pm\varepsilon_3\pm\varepsilon_4)$ 
is indicated by $1\pm2\pm3\pm4$. Thus the simple roots are 
$\alpha_1=2-3$, $\alpha_2=3-4$, $\alpha_3=4$, 
$\alpha_4=1-2-3-4$.
If $q$ is odd, then the possible representatives are  $x_i$, $i\in \I_{25}$, for $p\neq 3$, 
with two additional representatives $x_{26},x_{27}$ when $p=3$.

\begin{lema}\label{lem:f4} If $q$ is odd, then $\Oc$ is not kthulhu.
\end{lema}
\pf  
A  direct verification shows that all representatives for $i\geq 7$ enjoy the  $\alpha\beta$-property with $(\alpha, \beta)\neq 0$; we list
in Table~\ref{tab:f4-odd-char} the roots $\alpha$ and $\beta$ for each representative.
By Proposition \ref{prop:alfabetachevstein}, $\Oc$ is of type D.
\begin{table}[ht]
\caption{$\Oc_{x_i}$ with the $\alpha\beta$-property}\label{tab:f4-odd-char}
\begin{center}
\begin{tabular}{|c|c|c|}
\hline
$i$ &  $\alpha$ & $\beta$ \\
\hline
7&2&1-2+3+4\\
8&1-2&2-3\\
9,10&2-3&3-4\\
11&1+2-3-4&1-2+3+4\\
12,13&4&1-4\\
14,15&2-4&1-2\\
16&2-4&1-2+3+4\\
17,21,22&2-4&1-2-3+4\\
18&2&1-2\\
19,23,24,25,(26,27)&2-3&3-4\\
20&2&1-2-3-4\\
\hline
\end{tabular}
\end{center}
\end{table}

We next consider the representative $x_1$, that equals $x_\gamma(1)$ for a long root $\gamma$. 
By the discussion in \S \ref{subsubsec:steinberg-endo},
$\Oc_{x_1}$ contains an element in $\U_{\alpha_1}^F$, that lies 
in the subgroup of type $A_2$ generated by $\U_{\pm\alpha_1}, \U_{\pm\alpha_2}$.  Theorem \ref{thm:slsp} applies.

Finally, we deal with the $x_i$'s, $i\in \I_{2,6}$. 
Let $\Le_1$ be the standard Levi subgroup  (of type $B_3$)
generated by the root subgroups $\U_\gamma$, for $\gamma=\pm \alpha_1$, $\pm\alpha_2$,
$\pm\alpha_3$.
We claim that  all $x_i$, $i\in \I_{2,6}$, are conjugated to elements in $M_1$; 
then the result follows by  Proposition  \ref{prop:bn-qodd}. 
Indeed, $x_2, x_3$ lie in $\U^F_{1-2}\U^F_{1+2}$; thus
conjugating by $s_{1-3}s_{2-4}$, we get a representative in $\U^F_{3-4}\U^F_{3+4}$. 
Also $x_5,\,x_6$ lie in  $\U^F_{2-3}\U^F_{2+3}\U^F_4$, and $x_4=x_2(1)x_{3+4}(1)$, 
so they all lie in $M_1$.
\epf

\begin{lema}\label{lem:f4-char2} If $q$ is even, then $\Oc$ is not kthulhu.
\end{lema}

\pf The representative $x_1$, respectively $x_2$, is equal to $x_\gamma(1)$ for a short, respectively long, 
root $\gamma$. 
By the discussion in \S \ref{subsubsec:steinberg-endo}, $\Oc_{x_1}$ intersects $\U_{\alpha_3}^F$ and
$\Oc_{x_2}$ intersects $\U_{\alpha_1}^F$. 
Let $M = \langle \U^F_{\pm\alpha_3}, \U^F_{\pm\alpha_4}\rangle$ 
and $M'  = \langle \U_{\pm\alpha_1}^F, \U_{\pm\alpha_2}^F\rangle$, both  of type $A_2$. 
Then $\Oc_{x_1} \cap M$, respectively $\Oc_{x_2} \cap M'$, is a unipotent class corresponding to the partition 
$(2,1)$ in $M$, respectively $M'$.  By Theorem \ref{thm:slsp}, these classes are not kthulhu.

We consider now the classes labelled by $i \in \I_{20, 34}$.
Let  $\Pa_1$ be the standard parabolic subgroup with standard Levi $\Le_1$ as in the proof of 
Lemma \ref{lem:f4}.  Set $y_i=\pi_{1}(x_i)$. 
Then the class $\Oc_{y_i}^{M_1}$ satisfies the $\alpha\beta$-property; we list in 
Table~\ref{tab:f4char2-proj} the roots $\alpha$ and $\beta$ for each representative. 
Since $\Phi_{\Pi_1}$ is of type $B_3$, the group $[\Le_1,\Le_1]$ is isogenous to $\Sp_{6}(\kk)$. 
By Remark \ref{obs:sp-not-involution}, $\Oc_{y_i}^{M_1}$ is not kthulhu, hence neither is $\Oc$. 

\begin{table}[ht]
\caption{$\Oc_{y_i}^{M_1}$ with the $\alpha\beta$-property.}\label{tab:f4char2-proj}
\begin{center}
\begin{tabular}{|c|c|c|}
\hline
$i$ & $\alpha$ & $\beta$ \\
\hline
$i\in\I_{20,21}$&$\alpha_1$&$\alpha_2+\alpha_3$\\
$i\in\I_{22,23}$&$\alpha_3$&$\alpha_1+\alpha_2$\\
$i\in\I_{24,28}$&$\alpha_1+\alpha_2$&$\alpha_2+2\alpha_3$\\
$i\in\I_{29,34}$&$\alpha_1$&$\alpha_2$\\
\hline
\end{tabular}
\end{center}
\end{table}

We consider now the classes labelled by $i \in \I' = \{3,4, 7,8,12\} \cup \I_{14,19}$.
Let $\Pa_2$ be the standard parabolic subgroup  with standard Levi $\Le_2$ (of type $C_3$) associated with 
$\Pi_2=\{\alpha_2,\alpha_3,\alpha_4\}$; here $\Phi^+_{\Pi_2}$ consists of the roots $1-2$, $3$, $4$, $3\pm4$, 
$1-2\pm3\pm4$. 
Let $\beta_1=\alpha_4$, $\beta_2=\alpha_3$, $\beta_3=\alpha_2$ be the simple roots of $\Phi^+_{\Pi_2}$. 
Set $z_i=\pi_{2}(x_i)$.
Now $\Oc_{z_i}^{M_2}$ is a unipotent class in $\Sp_6(q)$. Let $\I'' = \I' -\{3,4\}$.
 Table~\ref{tab:f4char2-proj-sympl} lists the index  $i \in \I''$,  the support of $z_i$ 
 and the  partition associated to $\Oc_{z_i}^{M_2}$, 
 obtained from the Jordan form of $z_i$ in $\Sp_6(\kk)$. 
 Since the partition is always different from $(2, 1^4)$, 
 the label of the class in $\Sp_6(q)$ is never $W(1)\oplus V(2)$, 
 whence $\Oc_{z_i}^{M_2}$ is not kthulhu by Theorem \ref{thm:slsp}.
 The remaining classes in $\I'$ are represented by $x_3=x_1(1)x_{1+2}(1)$ and $x_4=x_{2+3}(1)x_1(1)$.
Let  $x=(\dot{s}_{1-3}\dot{s}_{2-4})\trid x_3\in\Oc_{x_3}^{\Gb}$ and 
$y=(\dot{s}_{2-3}\dot{s}_{1-2}\dot{s}_{3})\trid x_4\in \Oc_{x_4}^{\Gb}$.
Then $x\in \U_3\U_{3+4}$, so $x\in \U_{\beta_2+\beta_3}^F\U^F_{2\beta_2+\beta_3}\subset M_2$, 
$y\in \U_{1-2}\U_3$, so $y\in \U_{2\beta_1+2\beta_2+\beta_3}^F\U_{\beta_2+\beta_3}^F\subset M_2$. 
The partition associated to $x$, respectively $y$, as unipotent element in $\Sp_{6}(q)$  is $(2,2,1,1)$, 
respectively $(2,2,2)$.
Hence, neither $\Oc_{x_3}^{\Gb}$ nor $\Oc_{x_4}^{\Gb}$ is  kthulhu by Theorem \ref{thm:slsp}.

\begin{table}[ht]
 
\caption{$\supp z_i$ and its partition}\label{tab:f4char2-proj-sympl}
\begin{center}
\begin{tabular}{|c|c|c|}
\hline
$i$ &   $\supp z_i$ & partition \\
\hline
7,8,12,15&$\beta_1+2\beta_2+\beta_3$&$(2,2,1,1)$\\
14&$2\beta_1+2\beta_2+\beta_3$, $2\beta_2+\beta_3$,&$(2,2,1,1)$\\
16&$2\beta_1+2\beta_2+\beta_3$, $\beta_1+2\beta_2+\beta_3$,&$(2,2,1,1)$\\
17,18&$2\beta_1+2\beta_2+\beta_3$, $\beta_1+\beta_2$,&$(2,2,1,1)$\\
19&$2\beta_1+2\beta_2+\beta_3$, $2\beta_2+\beta_3$, $\beta_1+\beta_2+\beta_3$&$(2,2,2)$\\
\hline
\end{tabular}
\end{center}
\end{table}

The $x_i$'s for $i \in \I'''= \{5,6,9,10,11,13\}$ lie in the subgroup $\Kb$ of type $B_4$ 
generated by the subgroups $\U_{\pm\alpha}$,
$\alpha\in\{1-2, 2-3, 3-4,4\}$. 
If $i\in \I'''$, $\Oc^{\Kb^F}_{x_i}$ has the $\alpha\beta$-property, 
see Table \ref{tab:f4char2}.
Since $\SO_{9}(\kk)$ is isogenous to $\Sp_{8}(\kk)$, 
Remark \ref{obs:sp-not-involution} applies.

\begin{table}[ht]
\caption{$\Oc^{\Kb^F}_{x_i}$ with the $\alpha\beta$-property, $i \in \I''= \{5,6,9,10,11,13\}$.}\label{tab:f4char2}
\begin{center}
\begin{tabular}{|c|c|c|}
\hline
$x_i$ &  $\alpha$ & $\beta$ \\
\hline
5,\,6 & 2+3 & 1-3\\
9,\,10,\,13& 2 & 1-2 \\
11& 1-4& 3+4\\
\hline
\end{tabular}
\end{center}
\end{table}

\epf

\subsection{Unipotent classes in  $G_2(q)$}\label{subsec:g2}

Here we deal with unipotent classes in $G_2(q)$, $q>2$. As for $F_4(q)$, we shall use 
explicit representatives of the classes, the parabolics being too small. 
The list of representatives can be found in \cite{chang} when $p>3$ and in \cite{eno} 
otherwise; see \eqref{eq:G2-reps-3}, \eqref{eq:G2-reps-p>3}, \eqref{eq:G2-reps-2}. 

\begin{lema}\label{lem:g2-qodd} If  $q$ is odd, then $\Oc$ is not kthulhu.
\end{lema}
\pf Assume first $p> 3$. By \cite[Theorems 3.1, 3.2, 3.9]{chang} every non-trivial class of $p$-elements in $\Gb$ is 
either regular or can be represented by an element of the following form, for suitable $a,b,c\in \kc^{\times}$:
\begin{align}\label{eq:G2-reps-3}
\begin{aligned}
&x_{\alpha_2}(1),&  &x_{\alpha_2}(1)x_{3\alpha_1+\alpha_2}(b), & 
&x_{\alpha_2}(1)x_{2\alpha_1+\alpha_2}(-1)x_{3\alpha_1+\alpha_2}(c),\\
&x_{\alpha_1+\alpha_2}(1), & &x_{\alpha_2}(1)x_{2\alpha_1+\alpha_2}(a).& &
\end{aligned}
\end{align}

The regular classes  are covered by Proposition \ref{prop:regular} \eqref{item:regular-chevalley}.
The elements $x_{\alpha_2}(1)$ and $x_{\alpha_2}(1)x_{3\alpha_1+\alpha_2}(b)$  
lie in the  subgroup of type $A_2$ generated by $\U^F_{\pm\alpha_2}$ and $\U^F_{\pm(3\alpha_2+\alpha_2)}$ 
and we apply Theorem \ref{thm:slsp}. 
The classes represented by  $x_{\alpha_2}(1)x_{2\alpha_1+\alpha_2}(-1)x_{3\alpha_1+\alpha_2}(c)$ enjoy 
the $\alpha\beta$-property, so we invoke Proposition \ref{prop:alfabetachevstein}.
We prove now that the class of $r=x_{\alpha_1+\alpha_2}(1)$ is of type D. We observe 
that there is an element $\sigma = \dot s_{\alpha_{2}} \in\Gb\cap N_{\G}(\T)$ such that 
$s := \sigma\trid r=x_{\alpha_1}(\xi)$, $\xi\in \kc^{\times}$. 
Then  $sr\neq rs$ by the Chevalley commutator formula \eqref{eq:chev}
and, as $rs,sr\in\U^F$ and $p$ is odd, 
we have $(rs)^2\neq(sr)^2$. In addition, $r,s\in \Pa_1^F$, for $\Pa_1$ the standard parabolic subgroup 
with Levi $\Le_1$ associated with $\alpha_1$. Since $r$ lies in the unipotent radical $\Vu_1$ of $\Pa_1$ 
and $s$ lies in $\Le_1$, we have $\Oc_s^{\langle r,s\rangle}\neq   \Oc_r^{\langle r,s\rangle}$.

Let $r=x_{\alpha_2}(1)x_{2\alpha_1+\alpha_2}(a)$; it lies in 
$\langle \U_{\pm\alpha_2}^F\rangle\times\langle \U^F_{\pm (2\alpha_1+\alpha_2)}\rangle$. 
We argue as in \S \ref{subsec:further-rem}.  As $q>3$, Lemmata \ref{lem:x1x2y1y2} and 
\ref{lema:ACG-2-10} apply whence $\Oc_r$ is of type D.   

\medskip
Assume now $p=3$. By \cite[6.4]{eno} the non-trivial classes of $p$-elements in $\Gb$ are 
either regular or are represented by an element of the following form:
\begin{align}\label{eq:G2-reps-p>3}
\begin{aligned}
&x_{3\alpha_1+2\alpha_2}(1),&  &x_{\alpha_1+\alpha_2}(1)x_{3\alpha_1+\alpha_2}(a), \\
&x_{2\alpha_1+\alpha_2}(1)x_{3\alpha_1+2\alpha_2}(1), & &x_{2\alpha_1+\alpha_2}(1),
\end{aligned}
\end{align}
for suitable $a\in \kc^{\times}$. The regular classes  are covered by Proposition 
\ref{prop:regular} \eqref{item:regular-chevalley}. 
The element $x_{3\alpha_1+2\alpha_2}(1)$ lies in the subgroup of type $A_2$ generated 
by $\U^F_{\pm\alpha_2}$ and $\U^F_{\pm(3\alpha_1+2\alpha_2)}$ and  Theorem \ref{thm:slsp} applies.

We show that if $r=x_{\alpha_1+\alpha_2}(1)x_{3\alpha_1+\alpha_2}(a) \in \Oc$, then it  is of type D.
Indeed, let 
 $s:=\dot{s}_{\alpha_2}\trid r\in\U_{\alpha_1}^F\U_{3\alpha_1+2\alpha_2}^F$. Then $sr\neq rs$; 
 since $sr,rs\in\U^F$, we have $(sr)^2\neq(rs)^2$. 
Moreover, $r,s\in \Pa_1^F$ with $s\in\Le_1$, $r\in \Vu_1$, with notation as for $p > 3$. Thus,  
$\Oc_s^{\langle r,s\rangle}\neq   \Oc_r^{\langle r,s\rangle}$ and $\Oc$ is of type D.

Assume that $u=x_{2\alpha_1+\alpha_2}(1)x_{3\alpha_1+2\alpha_2}(1) \in \Oc$. 
Conjugating by suitable elements in $N_{\Gb}(\T)$ 
we find $r\in\Oc\cap\U_{\alpha_1}\U_{3\alpha_1+\alpha_2}\subset \Pa_1$, $r\not\in\Vu_1$ and 
$s\in\Oc\cap \U_{\alpha_1+\alpha_2}\U_{\alpha_2}\subset \Vu_1$. 
By repeated use of \eqref{eq:chev}, we see that the coefficient of $x_{\alpha_1+\alpha_2}$ in 
$srs^{-1}$ is $\neq 0$, hence $rs\neq sr$, $(rs)^2\neq(sr)^2$ and  $\Oc$ is of type D.

Assume finally that  $u=x_{2\alpha_1+\alpha_2}(1) \in \Oc$. 
Let $r=\dot{s}_{\alpha_1}\trid u\in\Oc_u^{\Gb}\cap\U_{\alpha_1+\alpha_2}$ and 
$s=\dot{s}_{\alpha_1+\alpha_2}\trid u\in\Oc_u^{\Gb}\cap\U_{\alpha_1}$. 
Then $rs,sr\in\U$, $(rs)^2\neq(sr)^2$, and $\Oc_s^{\langle r,s\rangle}\neq   \Oc_r^{\langle r,s\rangle}$, 
as $s\in \Le_1$ and $r\in \Vu_1$, so $\Oc$ is of type D.
\epf

In order to deal with some unipotent classes in $G_2(4)$ we will need a precise version of 
\eqref{eq:chev} for all pairs of positive roots. We shall use the relations from \cite[II.2]{eno}, 
that we write for convenience. They hold in general for $q$ even, and we shall use them 
recalling that $a^3=1$ for every $a\in\F_4^\times$. 
\begin{align} 
\label{eq:commutation-relations1}& x_{\alpha_1}(a)x_{\alpha_2}(b)  =
x_{\alpha_2}(b)x_{\alpha_1}(a)x_{\alpha_1+\alpha_2}(ab)x_{2\alpha_1+\alpha_2}(a^2b)x_{3\alpha_1+\alpha_2}(a^3b)\\
\label{eq:commutation-relations2}& x_{\alpha_1}(a)x_{\alpha_1+\alpha_2}(b) =
x_{\alpha_1+\alpha_2}(b)x_{\alpha_1}(a)x_{3\alpha_1+\alpha_2}(a^2b)x_{3\alpha_1+2\alpha_2}(ab^2)\\
\label{eq:commutation-relations3}& x_{\alpha_1}(a)x_{2\alpha_1+\alpha_2}(b)  =
x_{2\alpha_1+\alpha_2}(b) x_{\alpha_1}(a)x_{3\alpha_1+\alpha_2}(ab)\\
\label{eq:commutation-relations4}& x_{\alpha_2}(a)x_{3\alpha_1+\alpha_2}(b)  =
x_{3\alpha_1+\alpha_2}(b)x_{\alpha_2}(a)x_{3\alpha_1+2\alpha_2}(ab)\\
\label{eq:commutation-relations5}& x_{\alpha_1+\alpha_2}(a)x_{2\alpha_1+\alpha_2}(b)  =
x_{2\alpha_1+\alpha_2}(b)x_{\alpha_1+\alpha_2}(a)x_{3\alpha_1+2\alpha_2}(ab)
\end{align}
For all  other pairs of positive roots the corresponding subgroups commute.

\begin{lema}\label{lem:g2-char2} If $q > 2$ is even, then $\Oc$ is not kthulhu. 
\end{lema}
\pf  By \cite[2.6]{eno} all non-trivial classes of $2$-elements in $\Gb$ can be 
represented by an element of the following form, for suitable $a, b, c\in\kc$:
\begin{align}\label{eq:G2-reps-2}
\begin{aligned}
&x_{2\alpha_1+\alpha_2}(1),\quad x_{3\alpha_1+2\alpha_2}(1),& & 
x_{\alpha_1}(1)x_{\alpha_2}(1)x_{2\alpha_1+\alpha_2}(a), 
\\ & x_{\alpha_1+\alpha_2}(1)x_{2\alpha_1+\alpha_2}(1)x_{3\alpha_1+\alpha_2}(b), &
&x_{\alpha_2}(1) x_{2\alpha_1+\alpha_2}(1)x_{3\alpha_1+\alpha_2}(c).
\end{aligned}
\end{align}

Assume that $r=x_{2\alpha_1+\alpha_2}(1) \in \Oc$. It is enough to prove that $\Oc$ is of type C 
for $G_2(2)$, which is a non-simple subgroup of 
$G_2(q)$. We consider $\dot{s}_{\alpha_1+\alpha_2}\trid r=x_{\alpha_1}(1)\in\Oc$ and  
$s:=x_{-\alpha_1}(1)\trid x_{\alpha_1}(1)= \dot{s}_{\alpha_1}\in \Oc_r^{G_2(2)}$. 
Let $H:=\langle r,s, z=x_{\alpha_1+\alpha_2}(1)\rangle \le \Pa_1$ (the parabolic subgroup  
associated with $\alpha_1$), with $r\in \Vu_1$, $s\in\Le_1$. 
Hence, $\Oc_r^H\neq\Oc_s^H$. By a direct computation,
\begin{align*}
s\trid r &=x_{\alpha_1+\alpha_2}(1)=z\neq r,& z\trid r&=r x_{3\alpha_1+2\alpha_2}(1), \\ 
r\trid s &=szr,& z\trid(szr) &= s x_{3\alpha_1+2\alpha_2}(1). && 
\end{align*}
So  $H\leq \langle  \Oc_r^H,\Oc_s^H\rangle \leq H$;  $ \{r, z, z\trid r\}\subset\Oc_r^H$ and 
$\{s, szr, sx_{3\alpha_1+2\alpha_2}(1)\}\subset\Oc_s^H$ 
hence $\Oc_r^{G_2(2)}$ is of type C by \cite[Lemma 2.8]{ACG-III}. 

Assume that $r=x_{3\alpha_1+2\alpha_2}(1) \in \Oc$.
Now $r \in \M = \langle\U_{\pm\alpha_2}, \U_{\pm(3\alpha_1+\alpha_2)} \rangle$, of type $A_2$.
Since $\Oc_r^{\M}$ has partition $(2,1)$, $\Oc$ is not kthulhu by Theorem \ref{thm:slsp}
 and  \cite[Theorem 24.15]{MT}.

Assume $q>4$. The classes represented by  the $x_{\alpha_1}(1)x_{\alpha_2}(1)x_{2\alpha_1+\alpha_2}(a)$ 
for $a\in\kc$ are regular, thus they are not kthulhu by Proposition \ref{prop:regular}. 
The classes of $x_{\alpha_1+\alpha_2}(1)x_{2\alpha_1+\alpha_2}(1)x_{3\alpha_1+\alpha_2}(b)$   and
$x_{\alpha_2}(1) x_{2\alpha_1+\alpha_2}(1)x_{3\alpha_1+\alpha_2}(c)$ enjoy the $\alpha\beta$-property.  
By \cite[Proposition 3.6]{ACG-II}, these classes are of type F.

Let now $q=4$ and let  $\zeta$ be a generator of $\F_4^\times$ so $\zeta^2+\zeta+1=0$ and $\zeta^3=1$. 
By \cite{eno} there are $2$ regular unipotent classes, one represented by $x_{\alpha_1}(1)x_{\alpha_2}(1)$ 
and the other  by 
$x_{\alpha_1}(1)x_{\alpha_2}(1)x_{2\alpha_1+\alpha_2}(\zeta)$. We shall apply Lemma \ref{lem:typeC} 
in order to show that these classes are of type C.  For this, we need the following formula which can be 
retrieved applying  \eqref{eq:commutation-relations1} and \eqref{eq:commutation-relations4}.

\begin{align}\label{eq:abab}
\begin{aligned}
x_{\alpha_1}(a)x_{\alpha_2}(b)  x_{\alpha_1}(c)x_{\alpha_2}(d)  &=x_{\alpha_1}(a+c)x_{\alpha_2}(b+d)\\
 \times x_{3\alpha_1+\alpha_2}(b)&x_{3\alpha_1+2\alpha_2}(bd)x_{2\alpha_1+\alpha_2}(c^2b)  
 x_{\alpha_1+\alpha_2}(bc),  
\end{aligned}
\end{align}
$a,b,c,d \in \kc$.
Let  $r=x_{\alpha_1}(1)x_{\alpha_2}(1)$, $t:=\alpha_1^\vee(\zeta)$, 
$s:=t\trid r=x_{\alpha_1}(\zeta^2)x_{\alpha_2}(1)\in \Oc_r^{\Gb}$. 
By direct computation using \eqref{eq:abab} we see that
\begin{align*}
r^2&=x_{3\alpha_1+\alpha_2}(1)x_{3\alpha_1+2\alpha_2}(1)x_{2\alpha_1+\alpha_2}(1)  x_{\alpha_1+\alpha_2}(1)  \\
s^2&=x_{3\alpha_1+\alpha_2}(1)x_{3\alpha_1+2\alpha_2}(1)
x_{2\alpha_1+\alpha_2}(\zeta)  x_{\alpha_1+\alpha_2}(\zeta^2).
\end{align*} 
Using \eqref{eq:commutation-relations5} and that $\xi^{2}\neq \xi$,
we see $r^2s^2\neq s^2r^2$,  hence $r^2s\neq sr^2$ and  $s^2r\neq rs^2$. 
In addition, $\langle r,s\rangle\subseteq\U^F$ and 
$\U^F\trid r\subset r \langle \U_\gamma ~|~ \gamma\in \Phi^+ - \Delta\rangle$ and 
$\U^F\trid s\subset s \langle \U_\gamma~|~\gamma\in \Phi^+ - \Delta\rangle$, so 
$\Oc_r^{\langle r,s\rangle}\neq  \Oc_s^{\langle r,s\rangle}$, whence $\Oc_r^{\Gb}$ is of type C.

Similarly, we consider now $r=x_{\alpha_1}(1)x_{\alpha_2}(1)
x_{2\alpha_1+\alpha_2}(\zeta)$, $t:=\alpha_1^\vee(\zeta)$ 
and $s:=t\trid r=x_{\alpha_1}(\zeta^2)x_{\alpha_2}(1)x_{2\alpha_1+\alpha_2}(\zeta^2)\in \Oc_r^{\Gb}$.
In this case 
\begin{align*}
r^2&=x_{3\alpha_1+\alpha_2}(\zeta^2)x_{3\alpha_1+2\alpha_2}(\zeta^2)
x_{2\alpha_1+\alpha_2}(1)  x_{\alpha_1+\alpha_2}(1)  \\
s^2&=x_{3\alpha_1+\alpha_2}(\zeta^2)x_{3\alpha_1+2\alpha_2}(\zeta^2)
x_{2\alpha_1+\alpha_2}(\zeta)  x_{\alpha_1+\alpha_2}(\zeta^2).
\end{align*} As above we verify that $r^2s\neq s^2r$ and  $s^2r\neq r^2s$ and that 
$\Oc_r^{\langle r,s\rangle}\neq  \Oc_s^{\langle r,s\rangle}$ so $\Oc_r^{\Gb}$ is of type C.

We assume now that 
$x:=x_{\alpha_1+\alpha_2}(1)x_{2\alpha_1+\alpha_2}(1)x_{3\alpha_1+\alpha_2}(b) \in \Oc$, with $b\neq0$. 
By \cite[Proposition 2.6, page 499]{eno} if $q=4$ we can take $b = \zeta$. 
We prove that this class is of type C.
Set $r_\alpha:=x_{\alpha}(1)x_{-\alpha}(1)x_{\alpha}(1) = \dot s_{\alpha}$, $\alpha\in\Phi^+$, 
see \cite[Lemma 19]{yale}.
The elements
\begin{align*}
s=r_{\alpha_1}\trid x&=x_{2\alpha_1+\alpha_2}(1)x_{\alpha_1+\alpha_2}(1)x_{\alpha_2}(\zeta)\\
&=x_{\alpha_1+\alpha_2}(1)x_{2\alpha_1+\alpha_2}(1)x_{\alpha_2}(\zeta) x_{3\alpha_1+2\alpha_2}(1),
\\
r&=r_{\alpha_2}r_{\alpha_1}\trid s=x_{\alpha_1}(1)
x_{2\alpha_1+\alpha_2}(1)x_{3\alpha_1+2\alpha_2}(\zeta)
\end{align*}
 belong to $\Oc$.
We claim that $\Oc_r^{\langle r,s\rangle}\neq \Oc_s^{\langle r,s\rangle}$. 
Indeed,  $r, \,s\in\Pa_1$ with $r\not\in\V_1$, $s\in\V_1$.
A direct calculation shows that $r^2=x_{3\alpha_1+\alpha_2}(1)$, 
\begin{align*}
&r\trid s=x_{\alpha_1+\alpha_2}(1+\zeta)x_{2\alpha_1+\alpha_2}(1+\zeta)
x_{\alpha_2}(\zeta)x_{3\alpha_1+\alpha_2}(\zeta),\\
&r^2\trid s=s x_{3\alpha_1+2\alpha_2}(\zeta),\\
&r^3\trid s=r\trid(r^2\trid s)=r\trid (s x_{3\alpha_1+2\alpha_2}(\zeta))=(r\trid s)x_{3\alpha_1+2\alpha_2}(\zeta),\\
&s\trid(r\trid s)=x_{\alpha_1+\alpha_2}(1+\zeta)
x_{2\alpha_1+\alpha_2}(1+\zeta)x_{\alpha_2}(\zeta)x_{3\alpha_1+\alpha_2}(\zeta)x_{3\alpha_1+2\alpha_2}(1+\zeta).
\end{align*}
We see that all these are distinct, and different from $s$,
by looking at the unique expression as a product of elements in root subgroups in the order:
\begin{align*}
\alpha_1<\alpha_1+\alpha_2<2\alpha_1+\alpha_2<\alpha_2<3\alpha_1+\alpha_2<3\alpha_1+2\alpha_2
\end{align*}
Hence, $|\Oc_r^{\langle r,s\rangle}| \geq5$ and $\Oc$ is of type C, by \cite[Lemma 2.8]{ACG-III}, 
with $H=\langle r,s\rangle$.

%\comgas{[Here I am lost. Why do you consider the class of 
%$x=x_{\alpha_1+\alpha_2}(1)x_{2\alpha_1+\alpha_2}(1)$?]}

Let  now $t_2:=\alpha_1^\vee(\zeta)\alpha_2^\vee(\zeta)$, 
$t_3:=\alpha_2^\vee(\zeta)$, $t_4:=\alpha_1^\vee(\zeta)\alpha_2^\vee(\zeta^2)$ and set
\begin{align*}
x&=x_{\alpha_1+\alpha_2}(1)x_{2\alpha_1+\alpha_2}(1),\\
x_1&=r_{\alpha_2}\trid x=x_{\alpha_1}(1)x_{2\alpha_1+\alpha_2}(1)\in \Oc_x;
\\
x_2&=t_2\trid x_1=x_{\alpha_1}(\zeta)x_{2\alpha_1+\alpha_2}(\zeta),\\
x_3 &=t_3\trid x_1=x_{\alpha_1}(\zeta^2)x_{2\alpha_1+\alpha_2}(1),\\
x_4 &=t_4\trid x_1=x_{\alpha_1}(1)x_{2\alpha_1+\alpha_2}(\zeta).
\end{align*}
 Let $Y_i = \U^F\trid x_i$, $i \in \I_4$. A direct computation shows that
\begin{align*}
&Y_1= \bigcup_{f,\ell\in \F_4} x_{\alpha_1}(1)x_{\alpha_1+\alpha_2}(\ell)
x_{2\alpha_1+\alpha_2}(\ell+1)x_{3\alpha_1+2\alpha_2}(f^2+f)
\U_{3\alpha_1+\alpha_2}^F,\\
&Y_2= \bigcup_{f,\ell\in \F_4} x_{\alpha_1}(\zeta)x_{\alpha_1+\alpha_2}(\ell\zeta )
x_{2\alpha_1+\alpha_2}(\ell\zeta^2+\zeta)x_{3\alpha_1+2\alpha_2}(f^2\zeta+f\zeta)\U_{3\alpha_1+\alpha_2}^F,\\
&Y_3= \bigcup_{f,\ell\in \F_4} x_{\alpha_1}(\zeta^2)x_{\alpha_1+\alpha_2}(\ell^2\zeta)
x_{2\alpha_1+\alpha_2}(\ell\zeta+1)x_{3\alpha_1+2\alpha_2}(f^2\zeta^2+f)\U_{3\alpha_1+\alpha_2}^F,\\
&Y_4= \bigcup_{f,\ell\in \F_4} x_{\alpha_1}(1)x_{\alpha_1+\alpha_2}(\ell)
x_{2\alpha_1+\alpha_2}(\ell+\zeta)x_{3\alpha_1+2\alpha_2}(f^2+f\zeta)\U_{3\alpha_1+\alpha_2}^F.
\end{align*}
The union $Y = \bigcup_{i\in \I_4} Y_i$  is disjoint and a subrack of $\Oc_x$.
We take 
\begin{align*}
&r_1=x_1,&\\
&r_2=x_{\alpha_1}(\zeta)x_{\alpha_1+\alpha_2}(\zeta)x_{2\alpha_1+\alpha_2}(1)\in\U^F\trid x_2,&(\ell=1,f=0),\\
&r_3:= x_{\alpha_1}(\zeta^2)x_{2\alpha_1+\alpha_2}(1)\in\U^F\trid x_3,&(\ell=f=0),\\
&r_4:=x_{\alpha_1}(1)x_{\alpha_1+\alpha_2}(1)x_{2\alpha_1+\alpha_2}(\zeta^2)\in\U^F\trid x_4,&(\ell=1,f=0).
\end{align*}
We claim that  
$x_{\alpha_1}(a)x_{\alpha_1+\alpha_2}(b)x_{2\alpha_1+\alpha_2}(c)$ 
and $x_{\alpha_1}(\widetilde{a})x_{\alpha_1+\alpha_2}(\widetilde{b})x_{2\alpha_1+\alpha_2}(\widetilde{c})$ 
do not commute, for $a,b,c,\widetilde{a},\widetilde{b},\widetilde{c}\in\kc$ 
such that $c\widetilde{a}+\widetilde{a}^2b\neq \widetilde{c}a+a^2\widetilde{b}$. This follows from the formula:
\begin{align*}
&x_{\alpha_1}(a)x_{\alpha_1+\alpha_2}(b)x_{2\alpha_1+\alpha_2}(c)
x_{\alpha_1}(\widetilde{a})x_{\alpha_1+\alpha_2}(\widetilde{b})x_{2\alpha_1+\alpha_2}(\widetilde{c})=\\
&x_{\alpha_1}(a+\widetilde{a})x_{\alpha_1+\alpha_2}(b+\widetilde{b})
x_{2\alpha_1+\alpha_2}(c+\widetilde{c})x_{3\alpha_1+\alpha_2}
(c\widetilde{a}+\widetilde{a}^2b)x_{3\alpha_1+2\alpha_2}(b^2\widetilde{a}+c\widetilde{b}).
\end{align*}
 Hence, $r_ir_j\neq r_jr_i$ for $i\neq j$, $i,j\in\I_4$ and the class $\Oc_x$ is of type F. 

By \cite{eno}, the remaining class can be represented by any of
\begin{align*}
r_1&=x_{\alpha_2}(1) x_{2\alpha_1+\alpha_2}(1)x_{3\alpha_1+\alpha_2}(\zeta), \\
r_2 &= r_{\alpha_1} \trid r_1 = x_{\alpha_2}(\zeta) x_{\alpha_1+\alpha_2}(1)
x_{3\alpha_1+\alpha_2}(1)x_{3\alpha_1+2\alpha_2}(\zeta).
\end{align*}
Let $t:=\alpha_1^\vee(\zeta)\alpha_2^\vee(\zeta^2)$ and
\begin{align*}
&x:=t\trid r_1=x_{\alpha_2}(\zeta) x_{2\alpha_1+\alpha_2}(\zeta)x_{3\alpha_1+\alpha_2}(\zeta^2),\\
&y:=t\trid r_2=x_{\alpha_2}(\zeta^2) x_{\alpha_1+\alpha_2}(\zeta)x_{3\alpha_1+\alpha_2}(\zeta)
x_{3\alpha_1+2\alpha_2}(1).
\end{align*}

It is easier now to work with a different ordering of the positive roots:
\begin{align*}\alpha_1<\alpha_2<2\alpha_1+\alpha_2<\alpha_1+\alpha_2<3\alpha_1+\alpha_2<3\alpha_1+2\alpha_2.
\end{align*}
Let $Y_i= \U^F\trid r_i$, $i\in \I_2$, $Y_3 = \U^F\trid x$, $Y_4 = \U^F\trid y$.
A direct computation shows that
\begin{align*}
&Y_1= \bigcup_{\ell\in \F_4} x_{\alpha_2}( 1)x_{2\alpha_1+\alpha_2}(1+\ell^2)
x_{\alpha_1+\alpha_2}(\ell)x_{3\alpha_1+\alpha_2}(\zeta+\ell^3+\ell)\U_{3\alpha_1+2\alpha_2}^F,
\\
&Y_2= \bigcup_{\ell\in \F_4} x_{\alpha_2}( \zeta)x_{2\alpha_1+\alpha_2}(\ell^2\zeta )
x_{\alpha_1+\alpha_2}(1+\ell\zeta)x_{3\alpha_1+\alpha_2}(\ell^2+\ell^3\zeta+1)\U_{3\alpha_1+2\alpha_2}^F,\\
&Y_3= \bigcup_{\ell\in \F_4} x_{\alpha_2}(\zeta)x_{2\alpha_1+\alpha_2}(\zeta+\ell^2\zeta)
x_{\alpha_1+\alpha_2}(\ell\zeta)x_{3\alpha_1+\alpha_2}(\ell^3\zeta+\ell\zeta+\zeta^2)\U_{3\alpha_1+2\alpha_2}^F,\\
&Y_4= \bigcup_{\ell\in \F_4} x_{\alpha_2}( \zeta^2)x_{2\alpha_1+\alpha_2}(\ell^2\zeta^2)
x_{\alpha_1+\alpha_2}(\ell\zeta^2+\zeta)x_{3\alpha_1+\alpha_2}(\ell^3\zeta^2+\ell^2\zeta+\zeta)
\U_{3\alpha_1+2\alpha_2}^F.
\end{align*}
The union $Y = \bigcup_{i\in \I_4} Y_i$  is disjoint and a subrack of $\Oc$.
We take 
\begin{align*}
&r_3:= x_{\alpha_2}(\zeta)x_{2\alpha_1+\alpha_2}(1)x_{\alpha_1+\alpha_2}(1)\in Y_3, &(\ell=\zeta^2),\\
&r_4:=x_{\alpha_2}( \zeta^2)x_{2\alpha_1+\alpha_2}(\zeta^2)
x_{\alpha_1+\alpha_2}(1)x_{3\alpha_1+\alpha_2}(\zeta^2)\in
Y_4, &(\ell=1).
\end{align*}
By looking at the coefficient of $x_{3\alpha_1+2\alpha_2}$ in 
the expression of each product, we verify that $r_i\trid r_j\neq r_j\trid r_i$ if $i\neq j$, hence $\Oc$ is of type F.
\epf

\section{Unipotent classes in Steinberg groups}\label{sec:steinberg}

In this Section we deal with unipotent classes in Steinberg groups, i.e. $\PSU_n(q)$, $n \ge 3$ ;
$\Pom_{2n}^-(q)$, $n \ge 4$; ${}^3D_4(q)$ and ${}^2E_6(q)$. 
In order to apply inductive arguments as in Section \ref{sec:chevalley},
we first need information about the unitary groups $\PSU_n(q)$, including the non-simple group $\PSU_3(2)$.

\subsection{Unipotent classes in unitary groups}

Here $\Gb=\PSU_n(q)$, $G=\SU_n(q)$, $n\ge 3$ and $\G = \SL_n(\kk)$, for $n\geq 2$. 
For a clearer visibility of the behaviour of the conjugacy classes, we use the 
language of matrices and partitions. Here we choose $\B$, $\U$, as the subgroups of upper triangular, 
respectively unipotent upper triangular, matrices.
We start by some notation and basic facts. 

\begin{itemize}[leftmargin=*]\renewcommand{\labelitemi}{$\diamond$}
 \item ${\Jf}_n =\left(\begin{smallmatrix}
&&1\\
&\iddots&\\
1&&
\end{smallmatrix}\right) = {\Jf}_n^{-1} \in \GL_n(\kk)$.

\smallbreak
 \item $\Fr_q$ is the Frobenius endomorphism of $\GL_n(\kk)$ 
raising all entries of the matrix to the $q$-th power.

\smallbreak
 \item   $F: \GL_n(\kk) \to \GL_n(\kk)$, $F(X)={\Jf}_n\;^t\!(\Fr_q (X))^{-1}{\Jf}_n$, $X\in\GL_n(\kk)$. 

 \smallbreak
  \item $\GU_n(q) = \GL_n(\kk)^F$, $\SU_n(q) = \SL_n(\kk)^F\leq \SL_n(q^2)$, \cite[21.14(2), 23.10(2)]{MT}. 
%  Also, $\SU_n(q)$ can be realized as a subgroup of $\SL_n(q^2)$ \cite{W}. 
 
\smallbreak  \item 
To every  unipotent class in $\SU_n(q)$ we assign the partition of $n$ corresponding to the class 
in $\GL_n(q)$ it is embedded into. 
  
\smallbreak  \item 
Every  unipotent class in $\GL_n(\kk)$ 
meets $\GU_n(q)$ in exactly one class, since $C_{\GL_n(\kk)}(x)$ is connected for every $x$
\cite[8.5]{Hu}, \cite[I.3.5]{sp-st}. In other words, every partition comes from a class in $\SU_n(q)$.
 
\smallbreak \item  Since $\SU_n(q)$ is normal in $\GU_n(q)$, 
\cite[Remark 2.1]{ACG-I} says that 
all  unipotent class in $\SU_n(q)$ with the same partition are isomorphic as racks.

\smallbreak \item  For $d\leq n$ with $d\equiv n\mod 2$ and $h = \frac{n-d}2$, 
we denote by $\M_d\le \G$ the subgroup of matrices 
$\left(\begin{smallmatrix}
\id_{h}\\
& A\\
&&\id_{h}
\end{smallmatrix}\right)$ with $A\in\SL_d(\kk)$. So  $\M_d^F\simeq\SU_{d}(q)$.

\smallbreak \item  For $c\leq\left[\frac{n}{2}\right]$ we denote by $\Ha_{2c}\le \G$ 
the subgroup of matrices $\left(\begin{smallmatrix}
A&&B\\
&\id_{n-2c}\\
C&&D
\end{smallmatrix}\right)$ for $\left(\begin{smallmatrix}
A&B\\
C&D
\end{smallmatrix}\right)\in \SL_{2c}(\kk)$. Then  $\Ha_{2c}^F\simeq\SU_{2c}(q)$.

\smallbreak \item  If $q$ is odd, then $G^{\Fr_q}=\SO_n(q)$. If $q$  and $n$ are even, then $G^{\Fr_q}=\Sp_n(q)$.
\end{itemize}

\medbreak
Here is the main result of this Subsection:

\begin{prop}\label{lem:sungral}
Let $\Oc \neq \{e\}$ be a unipotent class in $\Gb = \PSU_n(q)$ with partition $\lambda$, 
where $\lambda$ is different from $(2,1,\ldots)$ if $q$ is even.  Then  $\Oc$ is not kthulhu. 
\end{prop}

\pf First, we reduce our analysis to $G = \SU_n(q)$ by the isogeny argument \cite[Lemma 1.2]{ACG-I}. Thus,
from now on $\Oc$ is a unipotent class in $G$.
Second, we split the proof for $q$ odd in \S \ref{subsec:odd-unitary} and for $q$ even in 
\S \ref{subsec:even-unitary}. In each of these, 
we distinguish several cases according to the partition 
$\lambda = (\lambda_1,\ldots,\lambda_n)$ associated to $\Oc$.

\subsubsection{Proof of Proposition \ref{lem:sungral} when $q$ is odd}\label{subsec:odd-unitary}

\begin{lema}\label{prop:sun-odd} If $\lambda_1 \ge 3$, but $\lambda \neq (3,1)$ in $\SU_4(3)$, 
then $\Oc$ is not kthulhu.
\end{lema}

\pf
If $\lambda_1$ is even, we may find $u$ in $\Oc\cap (\Ha^F_{\lambda_1}\times \M^F_{n-\lambda_1})$ 
such that its component in $\Ha^F_{\lambda_1}$ is regular. If $\lambda_1$ and $n$ are odd, then we may find $u$  
in $\Oc\cap (\Ha^F_{n-\lambda_1}\times \M^F_{\lambda_1})$ whose component in $\M^F_{\lambda_1}$ is regular. 
In both cases Proposition \ref{prop:regular} applies.

It remains the case when $\lambda_1$ is odd and $n$ is even. Then  there is $i > 1$ such that $\lambda_i$ is odd, 
and $l=\lambda_1+\lambda_i\geq4$. We take $i$ minimal with this property.  Then, we may 
find $u$ in $\Oc\cap(\Ha^F_{l}\times \M^F_{n-l})$, whose component in $\Ha^F_{l}$ 
has partition $(\lambda_1,\lambda_i)$. 
We consider $H:=\Ha^F_{l}\cap \Ha^{Fr_q}_{l}\cong\SO_l(q)$. Since $(\lambda_1,\lambda_i)$ 
is an orthogonal partition, we
may assume $u\in H$. If $l\geq 6$ the class $\Oc_u^H$ is of type D by Theorem \ref{thm:slsp}
and Proposition  \ref{prop:bn-qodd}, since $\SO_6(\kk)$ is $\SL_4(\kk)$ up to isogeny.
%\comgas{[It is not Proposition \ref{prop:DE}, instead?]}

Let $l=4$ and $q>3$. Now $\SO_l(\kk)$ is $\SL_2(\kk) \times \SL_2(\kk)$ up to isogeny, the class 
$\Oc_u^H$ is isomorphic as a rack to 
the product $X\times X$ for $X$ the non-trivial unipotent class in $\SL_2(q)$. 
By Lemma \ref{lema:ACG-2-10}, $\Oc_u^H$ is of type D. 

If $q=3$ and $n>4$, then the partition either contains the sub-partition $(3,3)$ or $(3,1,1,1)$. 
Reducing to the subgroup $\M_6$, we look at the classes $(3,3)$ and $(3,1,1,1)$ in $\SU_6(q)$. 
Since $\SO_6(q)<\SU_6(q)$ and $\SO_6(\kk)$ is $\SL_4(\kk)$ up to isogeny, these racks contain a subrack 
isomorphic to a non-trivial unipotent class in $\SL_4(q)$. Then Theorem \ref{thm:slsp} applies.
\epf

\begin{lema}\label{lem:3-1}If  $\lambda =(3,1)$,  $G = \SU_4(3)$, then $\Oc$ is of type D.
\end{lema}
\pf Let $\zeta$ be a generator of $\F_9^{\times}$. We may assume that 
$r:=\left(\begin{smallmatrix}
1&\zeta&\zeta&1\\
0&1&0&-\zeta^3\\
0&0&1&-\zeta^3\\
0&0&0&1
\end{smallmatrix}\right) \in \Oc$. Let $t:=\left(\begin{smallmatrix}
2&\zeta^6&\zeta^2&\zeta\\
2&\zeta^5&0&0\\
0&\zeta^2&2&\zeta^5\\
0&\zeta^6&1&\zeta^7
\end{smallmatrix}\right)\in \SU_4(3)$ and let $s=t\trid r= \left(\begin{smallmatrix}
0&1&0&1\\
2&2&1&2\\
0&0&0&2\\
0&0&1&2
\end{smallmatrix}\right)\in\Oc$. A direct computation shows that $(rs)^2\neq(sr)^2$.
The subgroup $H=\langle r,s\rangle\subset \{\left(\begin{smallmatrix}A&B\\
0&D\end{smallmatrix}\right)~|~A,D\in\SL_2(9)\}$. If $s\in\Oc_r^H$, then
$\left(\begin{smallmatrix}
0&1\\
2&2\\
\end{smallmatrix}\right)$ and $\left(\begin{smallmatrix}
1&\zeta\\
0&1\\
\end{smallmatrix}\right)$  would be conjugate in $\SL_2(9)$. But 
$\left(\begin{smallmatrix}
1&0\\
2&1\\
\end{smallmatrix}\right)\trid \left(\begin{smallmatrix}
0&1\\
2&2\\
\end{smallmatrix}\right)=\left(\begin{smallmatrix}
1&1\\
0&1\\
\end{smallmatrix}\right)$ which is not conjugate to $\left(\begin{smallmatrix}
1&\zeta\\
0&1\\
\end{smallmatrix}\right)$ because $\zeta$ is not a square. Hence, $\Oc_r^H\neq\Oc_s^H$ and $\Oc$  is of type D.
\epf

By Lemmata \ref{prop:sun-odd} and \ref{lem:3-1}, there remain the partitions $(2^a,1^b)$, $a > 0$. 

\begin{lema} If  $q>3$ and $\lambda_1=2$, then $\Oc$ is of type D.
\end{lema}
\pf Assume that $n$ is odd. Then  the partition contains $(2,1)$ and we may find a 
representative whose component in $\M_3^F$ has partition $(2,1)$. It is therefore enough to 
prove the statement for 
 $G=\SU_3(q)$ and $\lambda = (2,1)$. Let 
$r=\left(\begin{smallmatrix}
1&0&a\\
&1&0\\
&&1
\end{smallmatrix}\right)\in\Oc$ with $a\in\F_{q^2}^\times$, $a^q=-a$.
As $\kc^\times = \{\xi^{q+1} \vert  \xi\in\F_{q^2}^\times\}$ and $q>3$, we may pick  
$\xi\in\F_{q^2}^\times$ such that $-a^2\xi^{q+1} \in \kc^{\times} - (\{2\} \cup (\kc^{\times})^2)$. 
Let $t\in G$ be the diagonal matrix $(\xi, \xi^{q-1},\xi^{-q})$, 
$\sigma=\left(\begin{smallmatrix}
0&0&1\\
0&-1&0\\
1&0&0
\end{smallmatrix}\right)\in G$ and
$$s:=(\sigma t)\trid r=\left(\begin{smallmatrix}
1\\
0&1\\
a\xi^{1+q}&0&1
\end{smallmatrix}\right)\in\Oc.$$
Since $2 \neq -a^2\xi^{q+1}$, $(rs)^2\neq(sr)^2$. 
Let $\eta\in\kk$ be such that $\eta^2=a^{-1}$. 
Conjugating by the diagonal matrix $(\eta,\eta^{-1})$ we have
$$H:=\langle r,s\rangle \simeq 
\big\langle\left(\begin{smallmatrix}
1& a\\
0&1\end{smallmatrix}\right), \,\left(\begin{smallmatrix}
1& 0\\
a\xi^{q+1}&1\end{smallmatrix}\right)\big\rangle\simeq  
\big\langle\big(\begin{smallmatrix}
1& 1\\
0&1\end{smallmatrix}\big),\,
\left(\begin{smallmatrix}
1& 0\\
a^2\xi^{q+1}&1\end{smallmatrix}\right)\big\rangle.$$
By \cite[Theorem 6.21, page 409]{suzuki}, $H\simeq \SL_2(q)$. Since $-a^2\xi^{q+1}$ 
is not a square, $\Oc_r^H\neq\Oc_s^H$. 
Thus $\Oc$ is of type D.

Assume that $n$ is even. Then  the partition contains either $(2,2)$ or $(2,1,1)$ 
and we may use $\M_4$ to reduce  to  $\lambda = (2,2)$ or $(2,1,1)$ in $G=\SU_4(q)$. 
If  $\lambda = (2,2)$, which is an orthogonal partition, then we may assume that the 
representative $u$ lies in $\SO_4(q)$, and $\Oc_u^{\SO_4(q)}\cong X\times X$, where $X$ is the
 non-trivial unipotent class in $\SL_2(q)$. Hence, it is of type D.

If  $\lambda = (2,1,1)$, then we take
$r=\left(\begin{smallmatrix}
1&0&0&a\\
&1&0&0\\
&&1&0\\
&&&1
\end{smallmatrix}\right)\in\Oc$ for $a\in\F_{q^2}^\times$ satisfying $a^q=-a$.
Let $t\in G$ be the diagonal matrix $(\xi,\xi^{-1} ,\xi^{q},\xi^{-q})$ for $\xi\in\F_{q^2}^\times$, 
such that $-a^2\xi^{q+1}\in\kc^\times$ is not a square in $\kc^\times$, and let
$\sigma=\left(\begin{smallmatrix}
0&0&0&1\\
0&\zeta&0&0\\
0&0&\zeta^{-q}&0\\
1&0&0&0
\end{smallmatrix}\right)\in G$ with $\zeta\in \F_{q^2}^\times$ such that $\zeta^{q}=-\zeta$. 
We consider 
$$s:=(\sigma t)\trid r=\left(\begin{smallmatrix}
1\\
0&1\\
0&0&1\\
a\xi^{1+q}&0&0&1
\end{smallmatrix}\right)\in\Oc$$ and we proceed as for $n$ odd.
\epf

\begin{lema}If $q=3$ and $\lambda_1=\lambda_2=2$, then $\Oc$ is of type D.
\end{lema}
\pf By reducing to the subgroup $\Ha_4$, it is enough to consider $\lambda = (2,2)$ in $G =\SU_4(3)$.
Let $\zeta$ be a generator of $\F_9^\times$. Let 
\begin{align*}
r = \left(\begin{smallmatrix}1&\zeta& 0&0\\
&1&0&0\\
&&1&-\zeta^3\\
&&&1\end{smallmatrix}\right) \in \Oc, \  
\sigma &=\left(\begin{smallmatrix}
\zeta^2&0&0&0\\
0&0&1&0\\
0&1&0&0\\
0&0&0&\zeta^{2}
\end{smallmatrix}\right)\in G, \
s=\sigma\trid r=\left(\begin{smallmatrix}1&0& \zeta^3&0\\
&1&0&-\zeta\\
&&1&0\\
&&&1\end{smallmatrix}\right).
\end{align*}
It is easy to see that $(rs)^2\neq (sr)^2$. 
In addition,  $\langle r,s \rangle\subset\U^F$ and $\Oc_r^{\U^F}\neq\Oc_s^{\U^F}$, whence the statement.
\epf

\begin{lema}\label{lem:SUq3} If $\lambda = (2,1,\ldots)$ in $G = \SU_n(3)$, $n\geq3$, then $\Oc$ is of type C.
\end{lema}
\pf Let $\F_9^\times=\langle \zeta\rangle$. Without loss of generality we may assume that  
$$r=\left(\begin{smallmatrix}
1&&\zeta^2\\
&\ddots\\
&&1\end{smallmatrix}\right)=\id_n+\zeta^2 e_{1,n} \in \Oc.$$
% and  $a\in \F_{9}^\times$ is such that $a^2=-1$.
We consider, for $n$ odd, respectively even, the following element of $\SU_n(3)$:
\begin{align*}\sigma:=\left(\begin{smallmatrix}
0&0&0&0&\zeta\\
0&\id_{\frac{n-3}{2}}&0&0&0\\
0&0&-\zeta^2&0&0\\
0&0&0&\id_{\frac{n-3}{2}}&0\\
\zeta^{-3}&0&0&0&0\end{smallmatrix}\right)&& \tau:=\left(\begin{smallmatrix}
0&0&0&0&0&\zeta\\
0&\id_{\frac{n-4}{2}}&0&0&0&0\\
0&0&\zeta&0&0&0\\
0&0&0&\zeta^{-3}&0&0\\
0&0&0&0&\id_{\frac{n-4}{2}}&0\\
\zeta^{-3}&0&0&0&0&0
\end{smallmatrix}\right).\end{align*}
Accordingly, we set $s:=\sigma\trid r$ or $s:=\tau\trid r$. 
In both cases, $s=\left(\begin{smallmatrix}
1&&\\
&\ddots\\
\zeta^{-2}&&1\end{smallmatrix}\right)=\id_n+\zeta^{-2}e_{n,1}$. Then  $rs\neq sr$. Let 
$H:=\langle r,s\rangle$. We have
$$H\simeq \langle \left(\begin{smallmatrix} 1&\zeta^2\\0&1\end{smallmatrix}\right),\,\left(\begin{smallmatrix} 1&0\\ \zeta^{2}&1\end{smallmatrix}\right)\rangle\simeq\SL_2(3).$$
Conjugation by $\diag(\zeta^{-1},\zeta)$ and \cite[Theorem 6.21, page 409]{suzuki} give 
$H\simeq \langle \left(\begin{smallmatrix} 1&1\\0&1\end{smallmatrix}\right),\,\left(\begin{smallmatrix} 1&0\\ 1&1\end{smallmatrix}\right)\rangle\simeq\SL_2(3)$,
so $\Oc_r^H\neq\Oc_s^H$. We conclude by \cite[Lemma 2.7]{ACG-III}.
\epf

\subsubsection{Proof of Proposition \ref{lem:sungral} when $q$ is even}\label{subsec:even-unitary}

\begin{lema}\label{lem:regularSU} If $\Oc$ is a regular unipotent class, then it is not kthulhu. 
\end{lema}
\pf  By Proposition \ref{prop:regular} it is enough to deal with regular unipotent classes in $\SU_3(q)$ for $q=2,8$, $\SU_4(q)$ for $q=2,4$ and $\SU_n(2)$ for $n\geq 5$.  

\begin{enumerate}[leftmargin=*, label={(\roman*)}]
\item \emph{Regular unipotent classes  in $\SU_3(2^{2h+1})$,  $h\in{\mathbb N}_0$,  are of type D.}  
\end{enumerate}
It suffices to prove the claim for $G=\SU_3(2)\leq \SU_3(2^{2h+1})$. Let $\zeta$ be a generator of $\F_4^\times$. 
Consider the class $\Oc$ represented by
$r=\left(\begin{smallmatrix}
1&1&\zeta\\
0&1&1\\
0&0&1\\
\end{smallmatrix}\right)$. Let $t=\left(\begin{smallmatrix}
0&0&1\\
0&1&\zeta^2\\
1&\zeta&\zeta
\end{smallmatrix}\right)\in \SU_3(2)$ and $s:=t\trid r=\left(\begin{smallmatrix}
1&0&0\\
1&1&0\\
\zeta^2&1&1\\
\end{smallmatrix}\right)\in\Oc$. By direct verification, $(rs)^2\neq(rs)^2$. A computation with GAP shows that $\Oc_s^{\langle r,s\rangle}\neq  \Oc_r^{\langle r,s\rangle}$.
\begin{enumerate}[leftmargin=*, label={(\roman*)}]\setcounter{enumi}{1}
\item \emph{Regular unipotent classes in $\SU_n(q)$, $n\geq4$ even, are not kthulhu.} 
\end{enumerate}
Indeed, by the Jordan form theory, $\Oc$ is represented by an element
of a regular class in $\Sp_{n}(q)=\SU_n(q)^{\Fr_q}$. We conclude invoking  Theorem \ref{thm:slsp}.
\begin{enumerate}[leftmargin=*, label={(\roman*)}]\setcounter{enumi}{2}
\item  \emph{Regular unipotent classes in $\SU_n(2)$, $n \ge 5$ odd, are not kthulhu.}
\end{enumerate}
 By projecting a representative in $\U^F$ to $\M_5$, we obtain a regular unipotent class in the latter. Hence, it is enough to assume $n=5$. 
Let $\Pa$ be the standard $F$-stable parabolic subgroup  of $\G$ with Levi factor $\Le$ 
corresponding to the simple roots $\alpha_{2},\,\alpha_{3}$, 
and let $\pi\colon \Pa^F\to \Le^F$ be the projection. Any  $u\in\Oc\cap \U^F$ lies in $P:=\Pa^F$ and $\pi(u)$
is regular in $\Le^F$, which is the subgroup of matrices of the form
$\left(\begin{smallmatrix}
\det A\\
        & A\\
        && \det A
  \end{smallmatrix}\right)$, where $A\in \GU_3(2)$, so $\Le^F\simeq \GU_3(2)$. 
By Proposition \ref{prop:regular}, $\Oc_{\pi(u)}^{\Le^F}$ is not kthulhu. Then
Lemma \ref{lem:flexibility} applies. 
\epf

Now we argue inductively starting from Lemma \ref{lem:regularSU}.

\begin{lema}\label{lem:big-blocks-SU}
If any of the following conditions holds, then $\Oc$ is not kthulhu. 
\begin{enumerate}[leftmargin=*, label=\rm{(\alph*)}]
\item\label{item:su1} There is $j$ such that $\lambda_j\ge 4$ is even.
\item\label{item:su2} $n$ is odd and there is $j$ such that $\lambda_j \geq3$ is odd.
\item\label{item:su3} There is $i$ such that $\lambda_i=\lambda_{i+1}>2$.
\item\label{item:su6} $n$ is even and there are  $i,\,j$ such that $\lambda_j>\lambda_{i}\geq3$ and are both odd. 
\end{enumerate}
\end{lema}
\pf \ref{item:su1}, \ref{item:su2}: Apply Lemma \ref{lem:regularSU}  either to  $\Ha_{\lambda_j}$ or to $\M_{\lambda_j}$.

\ref{item:su3}: The class with partition $(\lambda_i,\lambda_i)$ in $\Ha_{2\lambda_i}^F$ has  a representative in 
$(\Ha_{2\lambda_i}^F)^{\Fr_q}\simeq \Sp_{2\lambda_i}(q)$.  Its class has label 
$W(\lambda_i)$ or $V(\lambda_i)\oplus V(\lambda_i)$ hence Theorem \ref{thm:slsp} applies.

\ref{item:su6}: By considering the class with partition $(\lambda_j,\lambda_i)$ in $\M_{\lambda_i+\lambda_j}$
we may assume $n=\lambda_i+\lambda_j$. Let $d=\lambda_j-\lambda_i$ and
let  $\Pa$ be the parabolic subgroup with standard Levi subgroup of type 
$A_{\lambda_i-1}\times A_{\lambda_i-1}$ associated with 
\begin{align*}
\{\alpha_k\in\Delta~|~k\in\I_{\lambda_i-1}\cup\I_{\lambda_j+1,\lambda_j+\lambda_i-1}\}.
\end{align*} 
Then $M\simeq \SL_{\lambda_i}(q^2)$. 
By Lemma \ref{lem:inductive-parabolic} and Theorem \ref{thm:slsp} it is enough to show that  $\Oc\cap V=\emptyset$. Now, if  $u\in V$ then it is  of the form 
\begin{align*}
\begin{pmatrix}
\id_{\lambda_i} & A_1 & A_2 \\ 0 & B & A_3 \\ 0& 0 & \id_{\lambda_i}
\end{pmatrix},
\end{align*}
for some  upper-triangular $B\in\SU_{d}(q)$ and some matrices $A_i$, $i\in\I_3$. Hence,   
\begin{align}\label{eq:V}
\operatorname{rk}(u-\id)&={\rm rk}\begin{pmatrix}
A_1 & A_2 \\  B -\id_{d}& A_3 
\end{pmatrix}\le d+\lambda_i.
\end{align}
On the other hand, if $u\in \Oc$ by Jordan form theory we have 
\begin{align}\label{eq:O}
\operatorname{rk}(u-\id)&=\lambda_i+\lambda_j-2=2\lambda_i+d-2.
\end{align}
As $\lambda_i\geq3$, conditions \eqref{eq:V} and \eqref{eq:O} are not compatible, so $\Oc\cap V=\emptyset$.
\epf

\begin{lema}\label{lem:22}
If there is $i$ such that $\lambda_i=\lambda_{i+1}=2$, then $\Oc$ is not kthulhu. 
\end{lema}
\pf Looking at  $\Ha_{2\lambda_i}^F$, we reduce to  the partition $(2,2)$ in $\SU_4(q)$. 
Let $u$ be a representative of a unipotent class with label $V(2)\oplus V(2)$ in $\Sp_4(q) = \SU_4(q)^{\Fr_q} \leq \SU_4(q)$. By Jordan form theory, 
we may assume that $u\in \Oc$. By Theorem \ref{thm:slsp}, $\Oc_u^{\Sp_4(q)}$ is not kthulhu, whence the statement.
%If $q>2$, then we proceed as in the proof of Lemma \ref{lem:big-blocks-SU} \ref{item:su3}. For $q=2$,  
%we take  
%$r=\left(\begin{smallmatrix}
%1&1&0&0\\
%0&1&0&0\\
%0&0&1&1\\
%0&0&0&1
%\end{smallmatrix}\right) \in \Oc$. Let $\zeta$ be a generator of $\F_4^\times$. Let 
%$s:=t\trid r=\left(\begin{smallmatrix}
%1&0&0&0\\
%1&1&1&0\\
%0&0&1&0\\
%0&0&1&1
%\end{smallmatrix}\right)\in\Oc$, where
%$t=\left(\begin{smallmatrix}
%0&\zeta&0&1\\
%\zeta^2&0&1&\zeta^2\\
%0&1&0&0\\
%1&0&0&1
%\end{smallmatrix}\right)\in \SU_4(2)$. Then $(rs)^2\neq(rs)^2$. With  GAP we check that $\Oc_s^{\langle r,s\rangle}\neq  \Oc_r^{\langle r,s\rangle}$ and $\Oc$ is of type D.
\epf

By Lemmata \ref{lem:regularSU}, \ref{lem:big-blocks-SU} and \ref{lem:22} there remain  
the partitions: $(2,1^a)$ for all $n \ge 3$, and  
$(\lambda_1,2,1^a)$, $(\lambda_1,1^a)$ for $\lambda_1>1$ odd and $n$ even.

\begin{lema}\label{lem:SU-neven-lambda21} If $n$ is even and  
$\lambda = (\lambda_1,2,1^a)$,  where $1<\lambda_1$ is odd, then $\Oc$ is of type D.
\end{lema}
\pf
It is enough to consider $\lambda = (\lambda_1,2,1)$. We pick a representative
of $\Oc$ lying in $\Ha_2^F\times\M_{\lambda_1+1}^F\simeq\SL_2(q)\times \SU_{\lambda_1+1}(q)$.
Then $\Oc$ contains a subrack isomorphic to $X\times Y$ where $X \neq \{e\}$ is a unipotent class in $\SL_2(q)$
and $Y$ is a unipotent class with partition $(\lambda_1,1)$ in $\SU_{\lambda_1+1}(q)$. 
The latter is not a class of involutions
because $\lambda_1>2$. By  \cite[1.4(ii)]{TZ} and Remark \ref{rem:real} there are $y_1\neq y_2\in Y$ such that 
$y_1y_2=y_2y_1$. Now Lemmata \ref{lem:x1x2} and \ref{lema:ACG-2-10} apply. \epf

\begin{lema}\label{lem:SU-neven-lambda>3} If $n$ is even and  
		$\lambda = (\lambda_1,1,\ldots)$ 
for some $3<\lambda_1$ odd, then $\Oc$ is not kthulhu.
\end{lema}
\pf It is enough to deal with the partition $(\lambda_1,1)$. Set $d:=(\lambda_1+ 1)/2> 2$.
Let  $\Pa$ be the parabolic subgroup with standard Levi subgroup associated with $\Delta- \alpha_{d}$. 
Then $M\simeq \SL_{d}(q^2)$. 
We claim that  $\Oc\cap V=\emptyset$. Indeed, if $u\in V$ then it is  of the form 
\begin{align*}
\begin{pmatrix}
\id_{d} & A \\  0 & \id_{d}
\end{pmatrix},
\end{align*}
for some $A$, so ${\rm rk}(u-1)\leq d$. If, in addition, $u\in \Oc$, then  ${\rm rk}(u- \id) =\lambda_1-1$. This is impossible because $\lambda_1 > 3$.
Lemma \ref{lem:inductive-parabolic} and Theorem \ref{thm:slsp} apply. 
\epf

\begin{lema}\label{lem:SU-neven-3111}If  $n\geq6$ is even and 
$\lambda = (3,1,1,1,\ldots)$, then $\Oc$ is of type D.
\end{lema}
\pf It is enough to deal with the partition $(3,1,1,1)$ in $\SU_6(q)$. 
Let $x\in \F_{q^2}^\times -  \kc^\times$ and let $r=\left(\begin{smallmatrix}
1&x&0&0&1&x\\
&1&0&0&0&1\\
&&1&0&0&0\\                                                   
&&&1&0&0\\
&&&&1&x^q\\
&&&&&1
\end{smallmatrix}\right)\in \Oc$,
$$\sigma= \left(\begin{smallmatrix}0&0&1&0&0&0\\
1&0&0&0&0&0\\                                                   
0&1&0&0&0&0\\
0&0&0&0&1&0\\
0&0&0&0&0&1\\
0&0&0&1&0&0\\
\end{smallmatrix}\right)=s_1s_5s_2s_4\in \SU_6(q),\quad 
s=\sigma\trid r=\left(\begin{smallmatrix}1&0&0&0&0&0\\
&1&x&1&x&0\\
&&1&0&1&0\\                                                   
&&&1&x^q&0\\
&&&&1&0\\
&&&&&1
\end{smallmatrix}\right)\in \Oc.$$
Since $r,s \in \U$, the discussion in  \cite[3.1]{ACG-I}  shows that 
$\Oc_r^{\langle r,s\rangle}\neq \Oc_s^{\langle r,s\rangle}$. 
By looking at the $(1,5)$-entry, we see that  $(rs)^2\neq(sr)^2$ and $\Oc$ is of type D. \epf

\begin{lema}\label{lem:SU-31-q>4} If $q>4$ and
 $\lambda = (3,1)$ in $\SU_4(q)$, then $\Oc$ is of type F.
\end{lema}
\pf Let $x,y\in \F_{q^2}$ such that $x^qy+y^qx\not=0$ and $\zeta_i\in\kc^\times$, 
for $i \in \I_4$, satisfying $\zeta_i\neq\zeta_j$ if $i\neq j$. 
Let
\begin{align*}
r&=\left(\begin{smallmatrix}
1&x&y&xy^q\\
&1&0&y^q\\
&&1&x^q\\
&&&1\end{smallmatrix}\right), &
t_i &=\left(\begin{smallmatrix}
1& & & \\
&\zeta_i& &\\
&&\zeta_i^{-1}& \\
&&&1\end{smallmatrix}\right), &
r_i&:=t_i\trid r=\left(\begin{smallmatrix}
1&x\zeta_i^{-1}&y\zeta_i&xy^q\\
&1&0&y^q\zeta_i\\
&&1&x^q\zeta_i^{-1}\\
&&&1\end{smallmatrix}\right).
\end{align*} 
Then $r, r_i \in \Oc$, since $t_i \in\SU_4(q)$. Let
$H:=\langle r_1,r_2,r_3,r_4 \rangle\subset \U^F$. By the discussion in \cite[3.1]{ACG-I} that 
$\Oc_{r_i}^H\neq\Oc_{r_j}^H$ if $i\neq j$. Then $r_ir_j=r_jr_i$ if and only if 
$(x^qy+xy^q)(\zeta_i\zeta_j^{-1}+\zeta_i^{-1}\zeta_j)=0$, if and only if $\zeta_i=\zeta_j$.
\epf

\begin{lema}The unipotent classes of type $(3,1)$ in $\SU_4(2)$ are of type D. 
\end{lema}
\pf Let $\zeta$ be a generator of $\F_4^\times$. We may assume that 
$r=\left(\begin{smallmatrix}
1&1&\zeta &\zeta\\
0&1&0&\zeta^2\\
0&0&1&1\\
0&0&0&1\end{smallmatrix}\right) \in \Oc$.
Let $t=\left(\begin{smallmatrix}
0&0&1&0\\
0&0&0&1\\
1&0&0&0\\
0&1&1&0\\
\end{smallmatrix}\right)\in \SU_4(2)$ and $s:=t\trid r=\left(\begin{smallmatrix}
1&1&0&0\\
0&1&0&0\\
\zeta^2&\zeta&1&1\\
0&\zeta&0&1
\end{smallmatrix}\right)\in\Oc$. Then $(rs)^2\neq(rs)^2$. By GAP we see that 
$\Oc_s^{\langle r,s\rangle}\neq  \Oc_r^{\langle r,s\rangle}$ and the claim follows.
\epf

\begin{lema}The unipotent classes of type $(3,1)$ in $\SU_4(4)$  are of type D. 
\end{lema}
\pf Let $\zeta$ be a generator of $\F_{16}^\times$. We may assume that
$r=\left(\begin{smallmatrix}
1&1&\zeta &\zeta\\
0&1&0&\zeta^4\\
0&0&1&1\\
0&0&0&1\end{smallmatrix}\right) \in \Oc$.
Let $t=\left(\begin{smallmatrix}
\zeta^{11}&\zeta^2&\zeta^5&\zeta^{14}\\
\zeta^{11}&\zeta^2&\zeta^8&\zeta^{2}\\
0&\zeta^{14}&\zeta^9&\zeta^{11}\\
0&\zeta^{14}&\zeta^9&\zeta^{6}\\
\end{smallmatrix}\right)\in \SU_4(4)$ and $s:=t\trid r=\left(\begin{smallmatrix}
0&1&0&\zeta^{12}\\
1&0&\zeta^3&\zeta^{10}\\
0&0&0&1\\
0&0&1&0
\end{smallmatrix}\right)\in\Oc$. We check at once that $(rs)^2\neq(rs)^2$, and with GAP that $\Oc_s^{\langle r,s\rangle}\neq  \Oc_r^{\langle r,s\rangle}$.
\epf
The Proposition is now proved.
\epf

The remaining classes could not be reached with our methods.

\begin{lema}\label{lem:SU211} If  $\lambda = (2,1^a)$ in $\SU_n(q)$, then $\Oc$ is austere, hence kthulhu.
\end{lema}

\pf We show that any subrack generated by two elements is either abelian or indecomposable. 
Let $r,s \in \Oc$, $rs \overset{\star}{\neq} sr$. We may assume 
$r=\id_n+a e_{1,n}=x_{\beta}(a)$ where $\beta$ is the highest positive root in $\Phi$ 
and $a\in \F_{q}^\times$. Let $g\in G$ be such that $s=grg^{-1}$. 
By \cite[24.1]{MT} there are $u,v\in\U^F$, and $\sigma\in G\cap N(\T)$ such that $g=u\sigma v$. 
As $F(\sigma)=\sigma$, the coset $\overline{\sigma}=\sigma \T\in W$ lies in $W^F\simeq \s_n^F$ 
which is the centralizer of the permutation 
\begin{align*}
(1,n)(2,n-1)\cdots ([\frac{n}{2}], n+1-[\frac{n}{2}]);
\end{align*}
hence, either $\overline{\sigma}(\{1,n\})=\{1,n\}$ or $\overline{\sigma}(\{1,n\})\cap \{1,n\}=\emptyset$.
Since $r$ is central in $\U^F$,
$s=u\sigma r \sigma^{-1} u^{-1}=u x_{\overline{\sigma}(\beta)}(a')u^{-1}$ for some $a' \in \kc$. 
Since $ru= ur$ and $rv=vr$, $\star$ holds if and only if 
$r\neq \sigma r \sigma^{-1}$.
Thus, $\overline{\sigma}(1)=n$ and $\overline{\sigma}(n)=1$, so
$\sigma$ is of the form
\begin{align*}
\sigma&=\left(\begin{smallmatrix}
0&0&\xi\\
0&A&0\\
\xi^{-q}&0&0
\end{smallmatrix}\right),& &\text{where } A\in \GU_{n-2}(q), &\xi&\in\kc^\times, \ \xi^{q-1}=\det A.
\end{align*}
Then $\sigma r\sigma^{-1}=\id_n+a\xi^{-1-q} e_{n,1}$, so
$$H:= \langle r,s\rangle\simeq
\langle \left(\begin{smallmatrix} 1&a\\0&1\end{smallmatrix}\right),\,\left(\begin{smallmatrix} 1&0\\ 
\xi^{-q-1} a&1\end{smallmatrix}\right)
\rangle\subset \SL_2(q).$$
Since the non-trivial unipotent class in $\SL_2(q)$ is sober \cite[3.5]{ACG-I}, $\Oc_r^H=\Oc_s^H$.
\epf

\subsection{Unipotent classes in $\Pom_{2n}^-(q)$, $n \ge 4$}

In this subsection $\Gb=\Pom_{2n}^-(q)$, $n \ge 4$.
We shall use the knowledge of unipotent conjugacy classes in $\PSL_n(q)$ and $\PSU_n(q)$ 
and apply inductive arguments. 

Here $\G$ is assumed simply-connected.
The root system of $\G$ is of type $D_n$, and the Dynkin diagram automorphism $\vartheta$ 
interchanges $\alpha_{n-1}$ and $\alpha_n$;  it fixes the basis vectors $\varepsilon_{j}$ for $j\in \I_{n-1}$, 
and maps $\varepsilon_n$ to $-\varepsilon_n$. 
Here is the main result of this Subsection:

\begin{prop}\label{prop:bn-steinberg-gral}
Let $\Oc$ be a non-trivial unipotent class in $\Pom^-_{2n}(q)$ with $n\geq 4$. Then $\Oc$ is not kthulhu. 	
\end{prop}

We split the proof for $q$ odd in \S \ref{subsubsec:pom-odd} and for $q$ even in \S \ref{subsubsec:pom-even}.

\subsubsection{Proof of Proposition \ref{prop:bn-steinberg-gral} when $q$ is odd}\label{subsubsec:pom-odd}

\pf Let $\Pa_1$ and $\Pa_2$ be the standard $F$-stable parabolic subgroups  with $F$-stable Levi factors 
$\Le_1$ and $\Le_2$ associated respectively with $\Pi_1:=\Delta -  \{\alpha_{n-1},\alpha_n\}$ (of type $A_{n-2}$), 
 and   $\Pi_2:=\{\alpha_{n-2}, \alpha_{n-1},\alpha_n\}$ (of type $A_{3}$). Then
\begin{align*}
\Phi_{\Pi_{1}}^{+} &= \{\varepsilon_{i}-\varepsilon_j ~|~ i< j\in \I_{n-1}\}, \ 
\Phi_{\Pi_{2}}^{+} = \{\varepsilon_{i}\pm \varepsilon_j ~|~ i< j\in \I_{n-2, n}\},\\
\Psi_{\Pi_{1}} &= \Phi^{+} \smallsetminus \Phi_{\Pi_{1}}^{+} = 
\{\varepsilon_{i}+\varepsilon_j, \varepsilon_{k}-\varepsilon_n ~|~i< j\in \I_{n},\ k\in \I_{n-1}\},
\\ 
\Psi_{\Pi_{2}} &=  
\{\varepsilon_{i}\pm \varepsilon_j ~|~   i< j,\ i\in \I_{n-3}, j \in \I_n\}.
\end{align*}  
By Lemma \ref{lem:inductive-parabolic}, Theorem \ref{thm:slsp} and Proposition \ref{lem:sungral}, 
it is enough to show
that $\Oc\cap\U^F \not\subset V_1\cap V_{2}$. 
Assume that there is  $u\in\Oc\cap V_1\cap V_{2}$; then 
\begin{align*}
\supp u\subset \Psi_{\Pi_1}\cap \Psi_{\Pi_2} = 
\{\varepsilon_{i}+\varepsilon_j, \varepsilon_{i}-\varepsilon_n ~|~ i< j,\ i\in \I_{n-3}, j \in \I_n\}.
\end{align*}
We consider various different cases.
\begin{enumerate}[leftmargin=*, label={(\roman*)}]\setcounter{enumi}{0}
	\item\label{it:bn-steinberg-gral-uno} $\varepsilon_{i}-\varepsilon_n \in \supp u$ for some 	$i\in \I_{n-3}$.
\end{enumerate}

Then $s_{\varepsilon_i-\varepsilon_{n-2}}(\supp u) \subseteq \Phi^{+}$. 
Since $s_{\varepsilon_{i}-\varepsilon_{n-2}}\in W^F$, it has a representative 
$\dot{s}_{\varepsilon_{i}-\varepsilon_{n-2}}\in N_{\G^F}(\T)$; hence 
$\dot{s}_{\varepsilon_{i}-\varepsilon_{n-2}}\trid u \in \Oc\cap \U^F$ and 
$\varepsilon_{n-2}-\varepsilon_{n}\in \supp(\dot{s}_{\varepsilon_{i}-\varepsilon_{n-2}}\trid u)$.
Thus $\dot{s}_{\varepsilon_{i}-\varepsilon_{n-2}}\trid u\in\U^\F\cap\Oc -  V_2$.

\begin{enumerate}[leftmargin=*, label={(\roman*)}]\setcounter{enumi}{1}
	\item $\varepsilon_{i}-\varepsilon_n \not\in \supp u$ for all $i\in \I_{n-3}$.
\end{enumerate}

Then  there exist $k\in  \I_{n-3}$ and $j$ such that $\varepsilon_{k}+\varepsilon_j \in \supp u$.
Let $$\ell = \max\{j~|~\ \varepsilon_{k}+\varepsilon_j \in \supp u\text{ for some } k\}.$$ If 
$\ell = n$, then pick a representative $\sigma\in N_{\G^F}(\T)\cap \Le_2$ 
of $s_{\varepsilon_{n-1}-\varepsilon_n}s_{\varepsilon_{n-1}+\varepsilon_n} \in W^F$. Thus 
$\sigma\trid u \in \Oc\cap V_{2}$ and 
$\varepsilon_{k}-\varepsilon_{n}\in \supp(\sigma \trid u)$ for all $k$ such that 
$\varepsilon_{k}+\varepsilon_n \in \supp u$. Therefore, 
either $\supp(\sigma\trid u)\not \subset V_1$, and we are done, or 
$\supp(\sigma\trid u)\subset V_1$ and $\varepsilon_k-\varepsilon_n\in\supp(\sigma\trid u)$, and 
we fall in \ref{it:bn-steinberg-gral-uno}.

If $\ell =n-1$, then pick a representative $\sigma\in N_{\G^F}(\T)\cap \Le_2$ 
of $s_{\varepsilon_{n-2}+\varepsilon_{n-1}} \in W^F$.
As above, $\sigma\trid u\in \Oc\cap V_{2}$, and 
$\varepsilon_i-\varepsilon_{n-2}\in\supp(\sigma\trid u)\cap \Phi_{\Pi_1}$ for some $i<n-2$. 
That is, $\supp(\sigma\trid u)\not \subset V_1$.

Finally, if $\ell<n-1$, then we pick  a representative $\sigma \in N_{\G^F}(\T)$ of 
$s_{\varepsilon_\ell-\varepsilon_{n-1}}$. 
Then  $\supp \sigma\trid u \subset  V_1\cap V_2$,  and we fall in the case $\ell=n-1$.
\epf

\subsubsection{Proof of Proposition \ref{prop:bn-steinberg-gral}, $q$  even}\label{subsubsec:pom-even}
Here Lemma \ref{lem:inductive-parabolic} does not apply in its full strength 
because of the existence of kthulhu classes in  $\PSU_4(q)$, $q$ even, and in $\PSL_3(2)$.
We proceed by induction on $n$. The case $n=4$,  Lemma \ref{lem:basis} below, requires a special treatment. 
\begin{lema}\label{lem:basis} 
If $\Gb = \Pom^-_{8}(q)$ with  $q$ even, then $\Oc$ is not kthulhu. 
\end{lema}
\pf Let us consider the $F$-stable standard parabolic subgroups $\Pa_1$, $\Pa_2$ 
with standard Levi subgroups $\Le_1$ and $\Le_2$
associated with the sets $\Pi_1=\{\alpha_1,\alpha_2\}$ and $\Pi_2=\{\alpha_2,\alpha_3,\alpha_4\}$, 
respectively. Let $u\in\Oc\cap\U^F$. 
We analyse different situations, according to $\Delta\cap\supp u$. Recall that, $u$ 
being $F$-invariant, the simple root $\alpha_3\in\supp u$ 
if and only if $\alpha_4\in\supp u$.

\begin{enumerate}[leftmargin=*, label={(\roman*)}]\setcounter{enumi}{0}
	\item $\alpha_2,\alpha_3,\alpha_4\in\supp u$.
\end{enumerate}
The projection $\pi_2(u)\in L_2$ is regular, thus $\Oc_{\pi_2(u)}^{M_2}$  is isomorphic as a rack 
to a unipotent class in $\SU_4(q)$ of partition $(4)$ and Proposition \ref{lem:sungral} applies.  

\begin{enumerate}[leftmargin=*, label={(\roman*)}]\setcounter{enumi}{1}
	\item $\Delta\cap\supp u= \{\alpha_1,\alpha_3,\alpha_4\}$ or 
	$\Delta\cap\supp u=\{\alpha_3,\alpha_4\}$.
\end{enumerate}

Then $\Oc_{\pi_2(u)}^{M_2}$ has partition $(2,2)$ or $(3,1)$ and Proposition \ref{lem:sungral} applies. 

\begin{enumerate}[leftmargin=*, label={(\roman*)}]\setcounter{enumi}{2}
	\item\label{it:alfa1} $\Delta\cap\supp u=\{\alpha_1\}$ or $\Delta\cap\supp u=\{\alpha_2\}$.
\end{enumerate}

Here $\pi_1(u)\in L_1$ is not regular, 
hence $\Oc_{\pi_1(u)}^{M_1}$  is isomorphic as a rack 
to a unipotent class in $\SL_3(q)$ with partition $\neq (3)$; Theorem \ref{thm:slsp} applies.
 
\begin{enumerate}[leftmargin=*, label={(\roman*)}]\setcounter{enumi}{3}
	\item\label{it:alfa12} $\Delta\cap\supp u=\{\alpha_1,\alpha_2\}$: either 
	$\alpha_2+\alpha_3+\alpha_4\in\supp u$ or 
	not.
 \end{enumerate}
We may assume that $\alpha_2+\alpha_3\not\in\supp u$, by conjugating with a suitable 
element in $(\U_{\alpha_3}\U_{\alpha_4})^F$ and using \eqref{eq:chev}. 
If $\alpha_2+\alpha_3+\alpha_4\in\supp u$, then $\Oc_{\pi_2(u)}^{M_2} \simeq \Oc_{v}^{\SU_4(q)}$, 
where ${\rm rk}(v-\id)=2$ and $(v-\id)^2=0$,
which is not kthulhu since its partition  is $(2,2)$. 
 If $\alpha_2+\alpha_3+\alpha_4\not\in\supp u$, then pick a representative
 $\sigma\in N_{\G^F}(\T)$  of $s_3s_4 \in W$. Then $\sigma\trid u\in\Oc\cap\U^F$ and 
 $\Delta\cap\supp(\sigma\trid u)=\{\alpha_1\}$ so we reduce to \ref{it:alfa1}.

\begin{enumerate}[leftmargin=*, label={(\roman*)}]\setcounter{enumi}{4}
	\item\label{it:alfa23} $\Delta\cap\supp u=\emptyset$ and $\alpha_1+\alpha_2\in \supp u$ 
	or $\alpha_2+\alpha_3\in \supp u$.
\end{enumerate}
In the first case,  $\Oc_{\pi_1(u)}^{M_1}$ has 
type $(2,1)$, and Theorem \ref{thm:slsp} applies. In the second, 
also  $\alpha_2+\alpha_4\in \supp u$ 
and $\Oc_{\pi_2(u)}^{M_2}$ has type $(2,2)$.
Indeed, $\Oc_{\pi_2(u)}^{M_2} \simeq \Oc_{v}^{\SU_4(q)}$, where ${\rm rk}(v-\id)=2$ and $(v-\id)^2=0$.
We invoke Proposition \ref{lem:sungral}. 

\begin{enumerate}[leftmargin=*, label={(\roman*)}]\setcounter{enumi}{5}
	\item\label{it:alfa234} $(\Delta \cup \{\alpha_1+\alpha_2, \alpha_2+\alpha_3,\alpha_2+\alpha_4\})\cap\supp u
	=\emptyset$.
\end{enumerate}

Let $\dot{s}_i\in N_{\G^F}(\T)$  be a representative of $s_i$, $i \in \I_2$.
If $\alpha_1+\alpha_2+\alpha_3\in \supp u$, then also $\alpha_1+\alpha_2+\alpha_4\in \supp u$.  
Now  $\dot{s}_1\trid u\in\U^F\cap\Oc$, $\Delta\cap\supp(\dot{s}_1\trid u)=\emptyset$ and  
$\alpha_2+\alpha_3\in \supp(\dot{s}_1\trid u)$, so we fall in \ref{it:alfa23}. Let $\sigma$ be as in \ref{it:alfa12}.
If $\alpha_2+\alpha_3+\alpha_4\in \supp u$, 
then $\sigma\trid u\in\Oc\cap \U^F$ and $\alpha_2\in\supp(\sigma\trid u)$ and we are in case \ref{it:alfa1}. 

\begin{enumerate}[leftmargin=*, label={(\roman*)}]\setcounter{enumi}{6}
	\item $\supp u\subset\{\alpha_1+\alpha_2+\alpha_3+\alpha_4, \alpha_1+2\alpha_2+\alpha_3+\alpha_4\}$. 
\end{enumerate}

If $\alpha_1+\alpha_2+\alpha_3+\alpha_4\in\supp u$, then  $\dot{s}_1\trid u$ is as in 
case \ref{it:alfa234}; while if $\supp u=\{\alpha_1+2\alpha_2+\alpha_3+\alpha_4\}$, 
then $\supp (\dot{s}_2\trid u) =\{\alpha_1+\alpha_2+\alpha_3+\alpha_4\}$. 
 \epf

We now proceed with the recursive step and assume that all non-trivial unipotent classes 
in a twisted group with root system $D_{n-1}$ are not kthulhu.

\smallbreak
Let $\Pa_1$ and $\Pa_2$ be the standard parabolic subgroups  
with  $F$-stable standard Levi subgroups $\Le_1$ and $\Le_2$
associated with the sets $\Pi_1=\{\alpha_i ~|~ i\in\I_{n-2}\}$ and $\Pi_2=\{\alpha_i ~|~ i\in\I_{2,n}\}$, 
of type $A_{n-2}$ and $D_{n-1}$ respectively. 
 By Lemma \ref{lem:inductive-parabolic} 
in order to prove the inductive step, it is enough to show that no non-trivial unipotent class $\Oc$ 
in $\G^F$ satisfies $\Oc\cap \U^F\subset V_1\cap V_2$. 
As usual let 
\begin{align*}
\Phi_{\Pi_1}&=\{\varepsilon_i-\varepsilon_j ~|~ i<j\in \I_{n-1}\}, \qquad\qquad\quad
\Phi_{\Pi_2}=\{\varepsilon_i\pm\varepsilon_j ~|~ i<j\in \I_{2,n}\},\\
\Psi_{\Pi_1}&=\{\varepsilon_i-\varepsilon_n,\,\varepsilon_j+\varepsilon_k~|~i\in\I_{n-1}, j<k\in \I_{n}\}, \,
\Psi_{\Pi_2} =\{\varepsilon_1\pm\varepsilon_j ~|~ j\in \I_{2,n}\}.
\end{align*}
Let $u\in\Oc\cap V_1\cap V_2$. 
Then $\supp u\subset \Psi_{\Pi_1}\cap\Psi_{\Pi_2}=
\{\varepsilon_1-\varepsilon_n, \varepsilon_1+\varepsilon_j~|~j\in \I_{2,n}\}$. 
Let $\dot{s}_i\in N_{\G^F}(\T)$  be a representative of  $s_i \in W^F$, $i \in \I_2$.
If $\supp u\neq\{\varepsilon_1+\varepsilon_2\}$ then $\dot{s}_1\trid u\in\Oc\cap \U^F$, 
but $\dot{s}_1\trid u\not\in V_2$  (look at its  support). If, instead, $\supp u=\{\varepsilon_1+\varepsilon_2\}$
then 
$\dot{s}_1\dot{s}_2\trid u\in\Oc\cap \U^F\cap \U_{\varepsilon_2+\varepsilon_3}$, 
so $\dot{s}_1\dot{s}_2\trid u\not\in V_1$. 

\medbreak This finishes the proof for $q$ even and 
Proposition \ref{prop:bn-steinberg-gral} is now proved. \qed

\subsection{Unipotent classes in ${}^2E_6(q)$}

We deal now with the group ${}^2E_6(q)$. Here  the Dynkin diagram automorphism $\vartheta$ 
interchanges $\alpha_{1}$ with  $\alpha_6$ and $\alpha_{3}$ with  $\alpha_5$. 
Here is the main result of this Subsection:

\begin{prop}\label{prop:E6-steinberg}
Let	$\Oc\neq \{e\}$  be a unipotent class in ${}^2E_6(q)$. Then $\Oc$ is not kthulhu. 	
\end{prop}

We give the proof for $q$ odd in \S \ref{subsubsec:2e6-odd} and for $q$ even in \S \ref{subsubsec:2e6-even}.
Let $\Pa_1$ and $\Pa_2$ be the $F$-stable standard parabolic subgroups  
with standard Levi subgroups $\Le_1$ and $\Le_2$
associated with $\Pi_1= \Delta -  \{\alpha_{2}\}$ (of type $A_{5}$)
and $\Pi_2=\{\alpha_{2}, \alpha_3, \alpha_{4}, \alpha_5\}$ (of type $D_{4}$).
Then  $\Psi_{\Pi_{1}}$, respectively $\Psi_{\Pi_{2}}$, consists of all positive roots containing $\alpha_{2}$,
respectively  at least  one of  $\alpha_{1}$ and $\alpha_{6}$.

\subsubsection{Proof of Proposition \ref{prop:E6-steinberg}, $q$ odd}\label{subsubsec:2e6-odd}

Here Lemma \ref{lem:inductive-parabolic} (c) applies softly to the  $F$-stable parabolic subgroups.

\pf  
 By Lemma \ref{lem:inductive-parabolic}, Propositions \ref{lem:sungral} and  \ref{prop:bn-steinberg-gral},
it is enough to show that $\Oc\cap\U^F \not\subset V_1\cap V_{2}$. Let 
$\beta = \sum_{i = 1}^{4}\alpha_i$, $\gamma =\sum_{i = 1}^{6}\alpha_i$;
thus $\vartheta(\beta) = \alpha_2 + \alpha_4 + \alpha_5 + \alpha_6$.
Let $u\in\Oc\cap \U^F$ lying in $V_1\cap V_2$.
Then  
\begin{align*}
\supp u\subset \Psi_{\Pi_{1}} \cap \Psi_{\Pi_{2}} = \Psi(\beta) \cup \Psi(\vartheta(\beta)) 
= \Sigma \cup \vartheta(\Sigma) \cup \Psi(\gamma);
\end{align*}
here $\Sigma = \{\beta_j~|~ j \in \I_{0,3}\}$ and $\Psi(\gamma) = \{\gamma_j~|~j \in \I_{0,6}\}$, 
where
\begin{align*}
\beta_0&=\beta, &\beta_1 &=s_5\beta_0; &\beta_2&=s_4\beta_1; &\beta_3&=s_3\beta_2;\\
\gamma_0&=\gamma, &\gamma_1 &=s_4\gamma_0; &\gamma_2&=s_3\gamma_1; &\gamma_3&=s_5\gamma_1;\\
\gamma_4&=s_5\gamma_2=s_3\gamma_3; &\gamma_5&=s_{4}\gamma_{4};& \gamma_6&=s_2\gamma_5. &&
\end{align*}

Let $\dot{w} \in N_{\G^F}(\T)$ be a representative  of $w \in W^F$. 
If either $\beta_j$ or $\vartheta(\beta_j) \in\supp u$ for $j \in \I_{0,3}$, 
then  $\dot{w_j}\trid u\in\Oc\cap\U^F - V_1$,
where $w_0 = w_1 = s_2$, $w_2 = s_2s_4$, $w_3 = s_2s_4s_5s_3$.
Thus we may assume that $\supp u \subset \Psi(\gamma)$.

If $\gamma_0\in\supp u$, then $\dot{s}_2\trid u\in\Oc\cap\U^F$, 
$\gamma-\alpha_2\in\supp(\dot{s}_2\trid u) - \Psi_{\Pi_1}$. Now we argue inductively. 
Suppose that  $\gamma_i\in\supp u$ for some $i\in \I_{0,j-1}$ implies that  $\Oc$ is not kthulhu. 
Assume that $\gamma_i\not\in\supp u$ for $i\in\I_{0,j-1}$ and $\gamma_j\in\supp u$. 
We claim that there is $\omega_j\in W^F$ with $\dot{\omega}_j\trid u\in\Oc\cap \U^F$ 
and either $\supp(\dot{\omega}_j\trid u)\not\subset\Psi(\gamma)$ (a case settled above), 
or $\gamma_l\in\supp(\dot{\omega}_j\trid u)$ for some $l\in\I_{0,j-1}$, where the recursive hypothesis applies. 
The claim holds, taking $\omega_1=\omega_5=s_4$, $\omega_2=\omega_3=s_1s_6$, $\omega_4=s_3s_5$,  $w_6=s_2$.  
\epf

\subsubsection{Proof of Proposition \ref{prop:E6-steinberg}, $q$ even}\label{subsubsec:2e6-even}

Here, the use of Lemma 
\ref{lem:inductive-parabolic} is hampered by the presence of kthulhu classes in  $\PSU_6(q)$.

\pf  As we have shown in the odd case, \S \ref{subsubsec:2e6-odd}, there is $u\in\Oc\cap\U^F$ 
such that $u\not\in V_1\cap V_2$. 
If $u\not\in V_2$, then  the result follows from Proposition \ref{prop:bn-steinberg-gral}.  
Let us assume  that $u\in V_2 - V_1$. In particular, $\alpha_3$,
$\alpha_4$, $\alpha_5$, $\alpha_3+\alpha_4$, $\alpha_4+\alpha_5$, $\alpha_3+\alpha_4+\alpha_5\not\in\supp u$.
Then  $\Oc_{\pi_1(u)}^{M_1}$ is non-trivial.

If $\supp u\cap \Phi_{\Pi_1}\neq\{\alpha_1+\alpha_3+\alpha_4+\alpha_5+\alpha_6\}$, then 
$\Oc_{\pi_1(u)}^{M_1} \simeq \Oc_{v}^{\SU_6(q)}$ where ${\rm rk }(v-\id)=2$, hence its 
associated partition is not $(2,1,1,1,1)$. By Lemma \ref{lem:sungral}, $\Oc_v^{M_1}$ and $\Oc$ are not kthulhu. 

If $\supp u\cap \Phi_{\Pi_1} = \{\alpha_1+\alpha_3+\alpha_4+\alpha_5+\alpha_6\}$, 
then $\dot{w}\trid u\in\Oc\cap \U^F$
where $w = s_1s_6\in W^F$. But $\alpha_3+\alpha_4+\alpha_5 \in \supp (\dot{w}\trid u)$, 
hence $\dot{w}\trid u \notin V_2$ and we are done. 
\epf

\subsection{Unipotent classes in ${}^3D_4(q)$}

We deal now with triality; $F$ arises from the graph automorphism $\vartheta$ of order $3$  
determined  by $\vartheta(\alpha_1)=\alpha_3$. We assume that $\G = \Gsc$.
We fix and ordering of the $\vartheta$-orbits in $\Phi^{+}$.
Let
\begin{align*}
y_{\alpha}(\xi) &:=x_{\alpha}(\xi)x_{\vartheta\alpha}(\xi^q)x_{\vartheta^2\alpha}(\xi^{q^2}),& 
\alpha&\in \Phi, \vartheta(\alpha)\neq\alpha, &\xi&\in \F_{q^3}.
\end{align*} 
Every element in $\U^F$ can be uniquely written as a product of elements $y_\alpha(\xi)$,
$\vartheta\alpha\neq\alpha$, $\xi\in\F_{q^3}$, and $x_\beta(\zeta)$, $\vartheta\beta=\beta$,  $\zeta\in\kc$.
Let 
\begin{align}\label{eq:triality-upsilon}
\Upsilon = \langle x_{\pm\gamma}(\xi), y_{\pm\delta}(\xi)~|~\vartheta(\gamma)=\gamma, 
\vartheta(\delta)\neq\delta, \xi\in\kc^\times  \rangle \leq \G^F.
\end{align}
The generators in \eqref{eq:triality-upsilon} are the non-trivial elements 
in the root subgroups with respect to  $\T^F\cap \T^{\Fr_q}$. 
It is known that $\Upsilon \simeq G_2(q) \simeq  G^{\Fr_q}$. 
 
\begin{prop}\label{prop:triality-unipotent} Every  unipotent class $\Oc \neq \{e\}$  in  ${}^3D_4(q)$ is not kthulhu.
\end{prop}
\pf By the isogeny argument we  work in $G = \G_{sc}^F$  \cite[Lemma 1.2]{ACG-I}. 
We analyse different cases separately, according to $q$ being odd, even and $>2$, or $2$.

\begin{enumerate}[leftmargin=*, label={(\roman*)}]
	\item \emph{$q$ is odd.} 
\end{enumerate}

The list of representatives of the unipotent classes in ${}^3D_4(q)$ appears in \cite[Table 3.1]{Ge};
they all have one of the following forms:
\begin{align*}
&x_{\alpha_1+2\alpha_2+\alpha_3+\alpha_4}(1), & &x_{\alpha_2}(1)y_{\alpha_1+\alpha_2+\alpha_3}(-1),
& u=&x_{\alpha_2}(1)y_{\alpha_1+\alpha_2+\alpha_3}(\zeta), \\
&y_{\alpha_1+\alpha_2+\alpha_3}(1), & & y_{\alpha_1}(1)x_{\alpha_2}(1),
& r=&y_{\alpha_1}(1)y_{\alpha_1+\alpha_2}(a),
\end{align*}
where $\zeta \in \F_{q^3}$ is not a square and $a\in \F_{q^3} - \kc$.
So  all classes but those of $u$ and $r$ have a representative in $\Upsilon \simeq G_2(q)$, 
hence they are not kthulhu by Lemma \ref{lem:g2-qodd}. Now 
$u\in H = \langle\U_{\pm\alpha_2}^F, y_{\pm(\alpha_1+\alpha_2+\alpha_3)}(b) ~|~ b\in\F_{q^3}^\times\rangle$, 
which is isogeneous to 
$\SL_2(q)\times \SL_2(q^3)$. Since $\Oc_u^H$ is the product of two non-trivial racks and $q^3>3$, 
$\Oc_u^H$ is of type D by Lemmata \ref{lema:ACG-2-10} and \ref{lem:x1x2y1y2}.

Assume that $r \in \Oc$.
Let $\xi$ be a generator of $\F_{q^3}^{\times}$, 
\begin{align*}
\eta&=\xi^{q-1},& t&=\alpha_1^\vee(\eta)\alpha_3^\vee(\eta^q)\alpha_4^\vee(\eta^{q^2}),& s&=t\trid r
=y_{\alpha_1}(\eta^2)y_{\alpha_1+\alpha_2}(a\eta^2).
\end{align*}
By \cite[Table 3.2]{Ge}, for every $b,c\in \F_{q^3}^\times$ we have
\begin{align}\label{eq:a,a+b}
\begin{aligned}
y_{\alpha_1+\alpha_2}(b)y_{\alpha_1}(c)
&=y_{\alpha_1}(c)y_{\alpha_1+\alpha_2}(b)y_{\alpha_1+\alpha_2+\alpha_3}(bc^q+cb^q) \\
&\ \times x_{\alpha_1+\alpha_2+\alpha_3+\alpha_4}(-(bc^{q^2+q}+b^qc^{q^2+1}+b^{q^2}c^{q+1})) \\
&\ \times x_{\alpha_1+2\alpha_2+\alpha_3+\alpha_4}(-(cb^{q^2+q}+c^qb^{q^2+1}+c^{q^2}b^{q+1})).
\end{aligned}
\end{align}
Using \eqref{eq:a,a+b} we verify that the coefficient of $y_{\alpha_1+\alpha_2+\alpha_3}$ in the expression of $rs$,
respectively $sr$, 
equals $a\eta^{2q}+a^q\eta^2$, respectively $a^q\eta^{2q}+a\eta^2$. 
These coefficients are equal if and only if $(a^q-a)(\eta^{2q}-\eta^2)=0$.
As $\eta^{2(q-1)}\neq1$ and $a^q\neq a$, we have $rs\neq sr$, with $rs,sr\in \U^F$. 
Thus, $(sr)^2\neq(rs)^2$, as $q$ is odd.
Comparing the coefficients of $x_{\alpha_1}$ in the expressions of $r$ and $s$ as products of 
elements in root subgroups, we see that 
$$
\U^F\trid r \subset x_{\alpha_1}(1)\langle \U_{\beta}|\beta\in\Phi^+-\{\alpha_1\}\rangle,\,
\U^F\trid s \subset x_{\alpha_1}(\eta^2)\langle \U_{\beta}| \beta\in\Phi^+-\{\alpha_1\}\rangle.
$$
So
$\Oc_r^{\langle r,s\rangle }\neq  \Oc_s^{\langle r,s\rangle }$,
whence $\Oc_r$ is of type D.

\begin{enumerate}[leftmargin=*, label={(\roman*)}]\setcounter{enumi}{1}
	\item \emph{$q > 2$ is even.} 
\end{enumerate}

The list of representatives of the unipotent classes in $G$ appears in \cite{DM}, 
see \cite[Table A2]{Hi}. For suitable $ \zeta,\zeta'\in\kc$,
the representatives are of the form
\begin{align*}
u_1&=x_{\alpha_1+2\alpha_2+\alpha_3+\alpha_4}(1),& 
u_2&=x_{\alpha_2}(1)x_{\alpha_1+\alpha_2+\alpha_3+\alpha_4}(1),\\
u_3&=y_{\alpha_1+\alpha_2+\alpha_3}(1), &
u_4&=y_{\alpha_1+\alpha_2}(1)y_{\alpha_1+\alpha_2+\alpha_3}(1)x_{\alpha_1+\alpha_2+\alpha_3+\alpha_4}(\zeta), \\
u_5&=y_{\alpha_1}(1)x_{\alpha_2}(1),&
u_6&=y_{\alpha_1}(1)x_{\alpha_2}(1)y_{\alpha_1+\alpha_2+\alpha_3}(\zeta'),\\
u_7&= y_{\alpha_1}(1)y_{\alpha_1+\alpha_2}(a),&  a &\in \F_{q^3} - \kc.
\end{align*}
All classes except those like $\Oc_{u_7}$ are represented by  $v \in \Upsilon \simeq G_2(q)$;
thus, these are not kthulhu by Lemma \ref{lem:g2-char2}. 
We deal with $\Oc_{u_7}$. Let  $\gamma_j = \sum_{i= 1}^{j} \alpha_i$ for shortness.
We use \eqref{eq:a,a+b} and the following relations from \cite{Ge}, cf. \cite{Hi}:
\begin{align*}
y_{\alpha_1}(b)y_{\gamma_3}(c) &= y_{\gamma_3}(c)y_{\alpha_1}(b)x_{\gamma_4}(c^qb+c^{q^2}b^q+cb^{q^2})\\
y_{\gamma_2}(b)y_{\gamma_3}(c) &= y_{\gamma_3}(c)y_{\gamma_2}(b)
x_{\alpha_1+2\alpha_2+\alpha_3+\alpha_4}(c^qb+c^{q^2}b^q+cb^{q^2}),\\
x_{\alpha_2}(d)x_{\gamma_4}(e) &= x_{\gamma_4}(e)x_{\alpha_2}(d)x_{\alpha_1+2\alpha_2+\alpha_3+\alpha_4}(de),\\
y_{\alpha_1}(b)x_{\alpha_2}(d) &= x_{\alpha_2}(d)y_{\alpha_1}(b)y_{\gamma_2}(bd)
y_{\gamma_3}(db^{q+1})x_{\gamma_4}(db^{q^2+q+1});
\end{align*}
here  $b,c\in \F_{q^3}^\times$ and $d, e \in \kc^\times$.
Let $\Ci \le \F_{q^3}^\times$ be the cyclic subgroup of order $q^2+q+1$ and $\Dc:=\Ci\cap\kc^\times$, a cyclic group of order $(q-1,3)$. Thus $|\Ci/ \Dc|  \ge 4$. Let  $\xi_i$, $i\in\I_{4}$, be representatives of $4$ distinct cosets in $\Ci/\Dc$ and let
\begin{align*}
t_i:&=\alpha_1(\xi_i)\alpha_3(\xi_i^q)\alpha_4(\xi_i^{q^2}),&
r_i:&=t_i\trid u_7 =y_{\alpha_1}(\xi_i^2)y_{\alpha_1+\alpha_2}(a\xi_i^{2})\in \Oc_r\cap \U^F.
\end{align*}
Since $\U^F\trid r_i\subset y_{\alpha_1}(\xi_i^2)\langle \U_\gamma~|~\gamma\in\Phi^+ - \Delta\rangle$, we have $\Oc_{r_i}^{\langle r_1,r_2,r_3,r_4\rangle}\neq \Oc_{r_j}^{\langle r_1,r_2,r_3,r_4\rangle}$ for $i\neq j$. In addition by \eqref{eq:a,a+b}  we see that
\begin{align*}
r_ir_j&\in y_{\alpha_1}(\xi_i^2+\xi_j^2)y_{\gamma_2}(a(\xi_i^2+\xi_j^2))y_{\gamma_3}(a\xi_i^2\xi_j^{2q}+a^q\xi_i^{2q}\xi_j^2)
\,\U^F_{\gamma_4}\U^F_{\alpha_1+2\alpha_2+\alpha_3+\alpha_4}
\end{align*}

The coefficients of $y_{\gamma_3}$ in the expressions of $r_ir_j$ and $r_jr_i$ are equal iff 
$(a+a^q)(\xi_i^2\xi_j^{2q}+\xi_j\xi_i^{2q})=0$, iff  $(\xi_i\xi_j^{-1})^{2(q-1)}=1$ (since $a\not\in\kc$),  iff $i=j$
by our choice of the $\xi_i$'s. Hence, $r_i\trid r_j\neq r_j$ for $i\neq j$ and $\Oc_{u_7}$ is of type F.

\begin{enumerate}[leftmargin=*, label={(\roman*)}]\setcounter{enumi}{2}
	\item \emph{$q = 2$.} 
\end{enumerate}

The description of the representatives is the same as in (ii) with $\zeta =0$ and $\zeta' = 1$, 
see \cite[\S 3]{Hi}, so that
\begin{align*}
u_4 &= y_{\alpha_1+\alpha_2}(1)y_{\alpha_1+\alpha_2+\alpha_3}(1), &
u_6&= y_{\alpha_1}(1)x_{\alpha_2}(1)y_{\alpha_1+\alpha_2+\alpha_3}(1).
\end{align*}
We do not have information on the unipotent classes of $G_2(2)$ yet, so we have to argue differently.
However, the argument for $u_7$ is exactly as for $q > 2$.
Now $u_1 \in \langle \U_{\pm\alpha_2}^F, \U_{\pm (\alpha_1+2\alpha_2+\alpha_3+\alpha_4)}^F \rangle$, 
a subgroup of type $A_2$, but it is not regular there. Hence
$\Oc_{u_1}$ is of not kthulhu by Theorem \ref{thm:slsp} and  \cite[Theorem 24.15]{MT}.

By \cite[Tables A.8]{Hi}, we have 
$r:=y_{\alpha_1}(1)y_{\gamma_3}(1)x_{\alpha_1+2\alpha_2+\alpha_3+\alpha_4}(1)\in \Oc_{u_2}$. 
Let $\xi\in \F_8^\times$ such that $\xi^3=\xi+1$. 
Then the roots  in $\F_8^\times$ of the polynomial $X^4+X^2+X$ are $\xi, \xi^2$ and $\xi^4$. 
Their inverses, together with $1$, are the roots of the polynomial
$X^4+X^2+X+1$. Let $\Pa_1$ be the parabolic subgroup with standard Levi subgroup 
associated with $\{\alpha_1,\alpha_3,\alpha_4\}$, and, 
 for $i\in\I_{4}$, let
\begin{align*}
&t_i:=\alpha_1^\vee(\xi^i)\alpha_3^\vee(\xi^{2i})\alpha_4^\vee(\xi^{4i}),\\
&r_i:=t_i\trid r=y_{\alpha_1}(\xi^{2i})y_{\gamma_3}(\xi^{6i})
x_{\alpha_1+2\alpha_2+\alpha_3+\alpha_4}(1)\in\Oc_{u_2},\mbox{ so }\\
&\U^F\trid r_a\subset y_{\alpha_1}(\xi^{2i})V_1.
\end{align*}
Then,  $\Oc_{r_i}^{\langle r_1,r_2,r_3,r_4\rangle}\neq \Oc_{r_j}^{\langle r_1,r_2,r_3,r_4\rangle}$ for $i\neq j$. 
In addition, 
\begin{align*}
r_ir_j&=y_{\alpha_1}(\xi^{2i}+\xi^{2j})y_{\gamma_3}(\xi^{6i}+\xi^{6j})x_{\gamma_4}(\xi^{4(j-i)}+\xi^{2(j-i)}+\xi^{j-i}).
\end{align*}
Let $i\neq j$. The coefficient of $x_{\gamma_4}$ in the expression of $r_ir_j$ is $0$ if and only if 
$\xi^{j-i}\in\{\xi,\xi^2,\xi^4\}$ if and only if the coefficient of $x_{\gamma_4}$ in the expression of $r_jr_i$ is $1$.  
Thus, $\Oc_{u_2}$ is of type F.

Let now $r_1=u_3$. Let $\sigma$ and $\tau$ in  $\Upsilon$ be representatives of $s_1s_3s_4, s_2 \in W^F$, 
respectively. Let $\Pa_2$ be the $F$-stable parabolic subgroup with standard Levi subgroup associated with 
$\alpha_2$.  We consider the following elements in $\Oc\cap V_2$:
\begin{align*}
r_2&=\sigma\trid r_1=y_{\alpha_1+\alpha_2}(1), \qquad
r_3=\tau\trid r_2=y_{\alpha_1}(1)\\
r_4&=x_{\alpha_2}(1)\trid r_3=y_{\alpha_1}(1)y_{\alpha_1+\alpha_2}(1)
y_{\alpha_1+\alpha_2+\alpha_3}(1)x_{\alpha_1+\alpha_2+\alpha_3+\alpha_4}(1).
\end{align*}
Let $Z= \left\langle \U_\gamma ~|~\gamma\in\Phi^+ - \{\alpha_1,\alpha_2,\alpha_1+\alpha_2\}\right\rangle$. Then  
\begin{align*}
V_2\trid r_1&\subset y_{\alpha_1+\alpha_2+\alpha_3}(1)Z,&
V_2\trid r_2&\subset y_{\alpha_1+\alpha_2}(1)Z, \\
V_2\trid r_3 &\subset y_{\alpha_1}(1) Z,&
V_2\trid r_4 &\subset y_{\alpha_1}(1)y_{\alpha_1+\alpha_2}(1) Z.
\end{align*}
Hence, the classes $\Oc_{r_i}^{\langle r_1,r_2,r_3,r_4\rangle}$ for $i\in\I_4$ are disjoint. 
A direct computation shows that $r_ir_j\neq r_jr_i$ for $i\neq j$, so $\Oc_{u_3}$ is of type F. 

We deal now with $u_4$. Let $\xi$, $\Pa_1$  and  $\Pa_2$ be as above and let
\begin{align*}
t_1&:=\alpha_1^\vee(\xi^3)\alpha_3^\vee(\xi^6)\alpha_4^\vee(\xi^5),&t_2&:=
\alpha_1^\vee(\xi)\alpha_3^\vee(\xi^2)\alpha_4^\vee(\xi^4),\\
r_1&:=t_1\trid u_4=y_{\gamma_2}(\xi^6)y_{\gamma_3}(\xi^4), &r_2&:=
x_{\alpha_2}(1)y_{\gamma_3}(1)y_{\gamma_4}(1),\\
r_3&:=y_{\alpha_1}(1)y_{\gamma_3}(1),&r_4&:=t_2\trid r_3=y_{\alpha_1}(\xi^2)y_{\gamma_3}(\xi^{-1}).
\end{align*}
Then, $r_i\in\Oc_{u_4}\cap\U^F$, \cite[Tables A.2, A.4, A.8, A.12]{Hi}. In addition,
\begin{align*}
\U^F\trid r_1&\subset V_1\cap V_2,&
\U^F\trid r_2&\subset x_{\alpha_2}(1) V_1\cap V_2, \\
\U^F\trid r_3 &\subset y_{\alpha_1}(1) V_1\cap V_2,&\U^F\trid r_4 &\subset y_{\alpha_1}(\xi^2)V_1\cap V_2.
\end{align*}
Hence, for $H=\langle r_i~|~i\in\I_4\rangle$ we have $\Oc_{r_i}^H\neq\Oc_{r_j}^H$ for $i,j\in\I_4$, with $i\neq j$. 
A direct computation shows that $r_ir_j\neq r_jr_i$, for $i\neq j$, so $\Oc_{u_4}$ is of type F.

Finally, we treat simultaneously the classes of $u_5$ and $u_6$, that are of the form
$x=y_{\alpha_1}(1)x_{\alpha_2}(1)y_{\alpha_1+\alpha_2+\alpha_3}(\epsilon)$ with
$\epsilon\in\{0,1\}$ respectively. Let $\Ci$ be as in the odd case and let $(\xi_i)_{i\in \I_{4}}$
be a family of distinct elements in $\Ci$. Set
\begin{align*}
t_i:&=\alpha_1^{\vee}(\xi_i)\alpha_3^{\vee}(\xi_i^q)\alpha_4^{\vee}(\xi_i^{q^2}),\\
r_i:&=t_i\trid x = y_{\alpha_1}(\xi_i^2)x_{\alpha_2}(1)
y_{\alpha_1+\alpha_2+\alpha_3}(\epsilon\xi_i^{1+q-q^2})\in \Oc_x\cap \U^F.
\end{align*}
Let $Q = \langle r_1,r_2,r_3,r_4\rangle$.
Since $\U^F\trid r_i\subset y_{\alpha_1}(\xi_i^2)x_{\alpha_2}(1)\langle 
\U_\gamma|~\gamma\in\Phi^+ - \Delta\rangle$, we have $\Oc_{r_i}^{Q}\neq \Oc_{r_j}^{Q}$ for $i\neq j$. 
The coefficient of $y_{\alpha_1+\alpha_2}$ in the expression of $r_ir_j$ equals 
$\xi_i^2$, hence $r_ir_j\neq r_jr_i$ for $i\neq j$. Hence $\Oc_{u_5}$ and $\Oc_{u_6}$ are  of type F.
\epf

%\]
\end{document}